\documentclass[english,twoside,11pt]{article}

\usepackage[german,english]{babel}
\usepackage{amsmath}
\usepackage{amssymb,latexsym}
\usepackage{moreverb}
\usepackage{url}
\usepackage{graphicx}
\usepackage{longtable}
\usepackage{lscape}
\usepackage{booktabs}
\usepackage{tabularx}
\usepackage{psfrag}
\usepackage{multicol}
\usepackage{paralist}

\newtheorem{deff}{Definition}
\newtheorem{lem}[deff]{Lemma}
\newtheorem{prop}[deff]{Proposition}
\newtheorem{thm}[deff]{Theorem}

\newtheorem{cor}[deff]{Corollary}

\makeatletter

\title{Equivelar and $d$-Covered Triangulations of Surfaces. II.\\ Cyclic Triangulations and Tessellations}

\author{\Large Frank H.~Lutz}

\date{.}

\begin{document}

\selectlanguage{english}

\maketitle

\begin{abstract}
With the $[0,1,2]$-family of cyclic triangulations we introduce a rich class of vertex-transitive triangulations of surfaces. 
In particular, there are infinite series of cyclic $q$-equivelar triangulations of orientable and non-orientable surfaces 
for every $q=3k$, $k\geq 2$, and every $q=3k+1$, $k\geq 3$. Series of cyclic tessellations of surfaces are derived 
from these triangulated series.
\end{abstract}

\section{Introduction}

A triangulation (as a finite simplicial complex) $K$ of a (closed) surface $M$ is \emph{equivelar} 
if all vertices of the triangulation have the same vertex-degree. An equivelar triangulation
with all vertices of degree~$q$ is called \emph{equivelar of type} $\{3,q\}$ or \emph{$q$-equivelar}.
An extensive survey of examples and basic properties of equivelar triangulations of surfaces 
is given in \cite{LutzSulankeTiwariUpadhyay2010pre}.

\enlargethispage*{.25mm}

The \emph{$f$-vector} (or \emph{face vector}) $f=(f_0,f_1,f_2)$ of a $q$-equivelar triangulation $K$
can be expressed as
\begin{equation}
f=(n,\frac{nq}{2},\frac{nq}{3}),
\end{equation}
where $n=f_0$, $f_1$, and $f_2$ are denoting the numbers of vertices, edges, 
and triangles of $K$, respectively.
For a $q$-equivelar triangulation of a surface $M$ of Euler characteristic $\chi(M)$,
\begin{equation}\label{eq:q}
q=6-\frac{6\chi(M)}{n}
\end{equation}
and
\begin{equation}\label{eq:bound_for_q}
q\leq\Bigl\lfloor\tfrac{1}{2}(5+\sqrt{49-24\chi (M)})\Bigl\rfloor.
\end{equation}
Since $q$ is a positive integer, it follows from (\ref{eq:q}) that for a
surface~$M$ of Euler characteristic $\chi(M)\neq 0$, the number of vertices $n$ 
has to be a divisor of $6|\chi(M)|$, and therefore $M$ has only finitely many equivelar triangulations.
Whereas (\ref{eq:bound_for_q}) gives an upper bound on~$q$ for equivelar triangulations,
Heawood's bound~\cite{Heawood1890} on the number of vertices $n$ from 1890 requires that 
\begin{equation}\label{eq:Heawood}
n\geq\Bigl\lceil\tfrac{1}{2}(7+\sqrt{49-24\chi (M)})\Bigl\rceil
\end{equation}
for any triangulation of a surface $M$ of Euler characteristic $\chi(M)$.
An integer triple $(\chi,q,n)$ is \emph{admissible} for a surface $M$
if $\chi=\chi(M)$ and if (\ref{eq:q}) and (\ref{eq:Heawood}) hold for $q$ and $n$, respectively.

Vertex-transitive triangulations are equivelar. However, based on enumerational results 
from \cite{KoehlerLutz2005pre}, \cite{Lutz1999}, and \cite{SulankeLutz2009},
it is remarked in \cite{LutzSulankeTiwariUpadhyay2010pre} that there are surfaces,
which do not have vertex-transitive triangulations for all their admissible triples.
In contrast, there is no case known where for an admissible triple $(\chi,q,n)$ of a surface $M$ 
with $\chi<0$ there are no equivelar triangulations. It is conjectured in \cite{LutzSulankeTiwariUpadhyay2010pre}
that respective equivelar triangulations always exist.

For each $q\geq 7$, McMullen, Schulz, and Wills \cite{McMullenSchulzWills1983} provide infinitely many 
examples of $q$-equivelar orientable triangulations. Moreover, they show 
that there are equivelar triangulations of orientable surfaces of genus $g$
with $O(g/{\log}\,g)$ vertices, which are realizable geometrically in ${\mathbb R}^3$.

Until recently, only two infinite series of vertex-transitive triangulations of surfaces
have been described explicitely in the literature, Ringel's cyclic $7{\rm mod}12$ series 
of neighborly orientable triangulations \cite{Ringel1961}, \cite{Ringel1974}
and Altshuler's cyclic series of triangulated tori \cite{Altshuler1971}.
In \cite{LutzSulankeTiwariUpadhyay2010pre}, a generalized Ringel series of cyclic orientable triangulations
$R_{k,n}$ with $k\geq 0$ and $n\geq 7+12k$ is given. The examples $R_{k,n}$ of this two-parameter family
have $n$ vertices, genus $g=kn+1$, and are equivelar with $q=6+12k$. The family $R_{k,n}$ 
contains for the boundary value $n=7+12k$ the subseries $R_{k,7+12k}$, which is Ringel's cyclic neighborly series, 
and for the boundary value $k=0$ the subseries $R_{0,n}$,
which is Altshuler's cyclic series of $6$-equivelar tori.

In this paper, a large class of cyclic series of triangulated surfaces is introduced in three steps.
In Section~\ref{sec:few_vertices}, a complete enumeration of all cyclic triangulations
of surfaces with up to $21$ vertices is given. From this data, various series of cyclic triangulations
are extracted. Section~\ref{sec:cyclic_series} presents respective series $A_{6}(n)$, $B_{9}(n)$, $C_{12}(n)$, 
$D_{12}(n)$, $E_{12}(n)$, $F_{12}(n)$, $G_{18}(n)$, $H_{18}(n)$, and $I_{18}(n)$.
Apart from the series $G_{18}(n)$, all these series belong to the $[0,1,2]$-family of cyclic triangulations,
which is defined in Section~\ref{sec:012}. Any cyclic triangulation of a surface is either $q$-equivelar 
with $q=3k$ or with $q=3k+1$ for some $k\geq 1$, as we will see in Corollary~\ref{cor:cyclic} of Section~\ref{sec:012}.
In Section~\ref{sec:012_series}, we construct two-parameter series
$S_{k,n}$, $T_{k,n}$, $\overline{T}_{k,n}$, $U_{k,n}$, $V_{k,n}$, $W_{k,n}$, $X_{k,n}$, $\overline{X}_{k,n}$, 
$Y_{k,n}$, and $Z_{k,n}$ of cyclic triangulations that belong to the $[0,1,2]$-family.
Moreover, we give a complete description of the case $q=7$. Together, this yields:

\bigskip


\noindent
\textbf{Main Theorem}\,\,
\emph{There are infinite series of cyclic $q$-equivelar triangulations of orientable and non-orientable surfaces
for each 
\begin{equation}
q = \left\{\begin{array}{lll}
                3k        & \mbox{\rm for} & k\geq 2,\\
                3k+1      & \mbox{\rm for} & k\geq 3.
               \end{array}\right.
\end{equation}
For $q=7$, there are exactly two cyclic $7$-equivelar triangulations, both of the
orientable surface of genus $g=2$ with $n=12$ vertices.
}

\bigskip

In Section~\ref{sec:beyond}, we describe the fundamental domains of cyclic triangulations 
that do not belong to the $[0,1,2]$-family. Section~\ref{sec:tessellations} is devoted to cyclic tessellations 
of surfaces that can be derived from cyclic triangulations. Particular classes of these tessellations 
have previously appeared in \cite{Brehm1990}, \cite{BrehmDattaNilakantan2002}, \cite{Datta2005}, and \cite{Jamet2001}.

\vspace{12mm}

\mbox{}

\pagebreak

\section{Cyclic triangulations of surfaces with few vertices}
\label{sec:few_vertices}

\begin{table}
\centering
\defaultaddspace=0.15em
\caption{Numbers of cyclic triangulations with up to $21$ vertices.}\label{tbl:cycl}
\begin{tabular*}{\linewidth}{@{}c@{\extracolsep{16pt}}c@{\extracolsep{16pt}}r@{\extracolsep{16pt}}r@{\extracolsep{\fill}}c@{\extracolsep{16pt}}c@{\extracolsep{16pt}}r@{\extracolsep{16pt}}r@{}}
\\\toprule
 \addlinespace
 \addlinespace
 \addlinespace
 \addlinespace
  Vertices  &  Orient.  & Genus  &  Types   &  Vertices  &  Orient.  & Genus  &  Types \\
\midrule
\\[-4mm]
 \addlinespace
 \addlinespace
 \addlinespace
 \addlinespace

      4    &        $+$ &      0 &       1  &     18    &        $+$ &      1 &       4  \\
      6    &        $+$ &      0 &       1  &           &            &      7 &       3  \\
      7    &        $+$ &      1 &       1  &           &            &     10 &       8  \\
      8    &        $+$ &      1 &       1  &           &            &     16 &       6  \\      9    &        $+$ &      1 &       1  &           &        $-$ &      2 &       1  \\

     10    &        $+$ &      1 &       1  &           &            &     14 &       7  \\
           &        $-$ &      2 &       1  &           &            &     20 &      22  \\
     11    &        $+$ &      1 &       1  &           &            &     32 &      30  \\
     12    &        $+$ &      1 &       4  &     19    &        $+$ &      1 &       3  \\
           &            &      2 &       2  &           &            &     20 &      32  \\

           &            &      5 &       3  &           &        $-$ &     21 &       6  \\
           &        $-$ &      8 &       5  &           &            &     40 &      46  \\
     13    &        $+$ &      1 &       2  &     20    &        $+$ &      1 &       6  \\
           &        $-$ &     15 &       2  &           &            &     16 &      76  \\
     14    &        $+$ &      1 &       2  &           &            &     21 &      78  \\

           &            &      8 &       2  &           &        $-$ &     12 &      10  \\
           &        $-$ &      2 &       1  &           &            &     22 &      36  \\
           &            &     16 &       5  &           &            &     32 &     192  \\
     15    &        $+$ &      1 &       4  &           &            &     42 &     162  \\
           &        $-$ &     12 &       5  &     21    &        $+$ &      1 &       6  \\

           &            &     17 &       5  &           &            &     22 &     250  \\
     16    &        $+$ &      1 &       3  &           &        $-$ &     16 &       5  \\
           &            &      5 &       5  &           &            &     23 &      11  \\
           &        $-$ &     18 &      11  &           &            &     37 &     176  \\
           &            &     26 &      32  &           &            &     44 &     287  \\

     17    &        $+$ &      1 &       2  & \\
           &        $-$ &     19 &       4  & \\
 \addlinespace

 \addlinespace
 \addlinespace
 \addlinespace
 \addlinespace
 \bottomrule
\end{tabular*}
\end{table}

A complete enumeration of all vertex-transitive triangulations of surfaces with up to $15$ vertices
is given in \cite{KoehlerLutz2005pre}, \cite{Lutz1999} by using the program \cite{Lutz_MANIFOLD_VT};
lists of all the vertex-transitive examples with up to $20$ vertices can be found online at \cite{Lutz_PAGE}.

We made use of the program \cite{Lutz_MANIFOLD_VT}
to search for triangulations of surfaces with vertex-transitive cyclic symmetry.
Table~\ref{tbl:cycl} lists the different types of cyclic examples with up to $21$ vertices.

\begin{thm}
Altogether, there are exactly $1570$ combinatorially distinct triangulated surfaces with up to $21$ vertices
that have a vertex-transitive cyclic symmetry. Of these examples, $508$ are orientable and $1062$ are non-orientable.
\end{thm}

\section{Some cyclic series with \mathversion{bold}$q=6,9,12,18$\mathversion{normal}}
\label{sec:cyclic_series}

Let the (action of the) cyclic group ${\mathbb Z}_n$ on the elements $0,1,\dots,n-1$ be generated by the permutation \mbox{$(0,1,\dots,n-1)$}.
For any cyclic triangulation of a surface, it suffices to give one generating triangle for every orbit of triangles.
Usually, the lexicographically smallest triangle of an orbit is chosen to represent the orbit.
For example, the triangle $[0,1,2]$ generates the orbit of triangles $[0,1,2]$, $[1,2,3]$, $[2,3,4]$, \dots, $[0,n-2,n-1]$, $[0,1,n-1]$. 
From the generating triangles of a cyclic triangulation we can easily read off the link of the vertex $0$ in the triangulation.
Moreover, we can determine $q$ and the $f$-vector of the triangulation from the number of distinct orbits and their sizes.
In the following, we present cyclic series $A_{6}(n)$, $B_{9}(n)$, $C_{12}(n)$, $D_{12}(n)$, $E_{12}(n)$, 
$F_{12}(n)$, $G_{18}(n)$, $H_{18}(n)$, and $I_{18}(n)$, which were found by analyzing the data of Section~\ref{sec:few_vertices}.
We defer a discussion of the orientability respectively non-orientability of the examples of the series to Section~\ref{subsec:or_vs_nonor}.

\subsection{The series \mathversion{bold}$A_{6}(n)$\mathversion{normal} of cyclic torus and Klein bottle triangulations}
\label{subsec:A_6_n}

\begin{figure}
\begin{center}
\small
\psfrag{0}{0}
\psfrag{1}{1}
\psfrag{2}{2}
\psfrag{3}{3}
\psfrag{n-2}{$n-2$}
\psfrag{n-1}{$n-1$}
\psfrag{n-4/2}{$\frac{n-4}{2}$}
\psfrag{n-2/2}{$\frac{n-2}{2}$}
\psfrag{n/2}{$\frac{n}{2}$}
\psfrag{n+2/2}{$\frac{n+2}{2}$}
\psfrag{n+4/2}{$\frac{n+4}{2}$}
\psfrag{n+6/2}{$\frac{n+6}{2}$}
\psfrag{n+8/2}{$\frac{n+8}{2}$}
\psfrag{n+10/2}{$\frac{n+10}{2}$}
\includegraphics[width=0.72\linewidth]{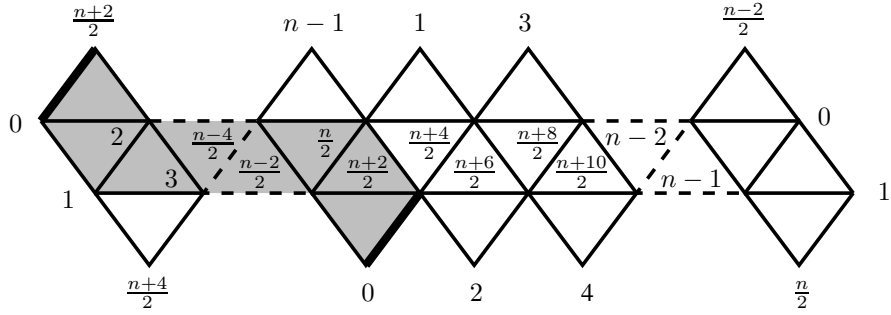}
\end{center}
\caption{The torus $A_{6}(8+4r)$.}
\label{fig:A_6_8_4r}
\end{figure}

\begin{figure}
\begin{center}
\small
\psfrag{0}{0}
\psfrag{1}{1}
\psfrag{2}{2}
\psfrag{3}{3}
\psfrag{n-2}{$n-2$}
\psfrag{n-1}{$n-1$}
\psfrag{n-6/2}{$\frac{n-6}{2}$}
\psfrag{n-4/2}{$\frac{n-4}{2}$}
\psfrag{n-2/2}{$\frac{n-2}{2}$}
\psfrag{n/2}{$\frac{n}{2}$}
\psfrag{n+2/2}{$\frac{n+2}{2}$}
\psfrag{n+4/2}{$\frac{n+4}{2}$}
\psfrag{n+6/2}{$\frac{n+6}{2}$}
\psfrag{n+8/2}{$\frac{n+8}{2}$}
\includegraphics[width=0.72\linewidth]{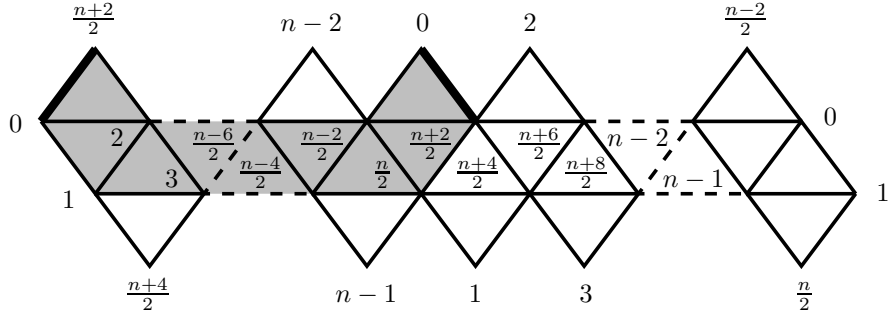}
\end{center}
\caption{The Klein bottle $A_{6}(10+4r)$.}
\label{fig:A_6_10_4r}
\end{figure}

Let $n=8+2m$ for $m\geq 0$. Then the two generating triangles
\begin{center}
\small
\begin{tabular}{c@{\hspace{4mm}}c@{\hspace{4mm}}c}
$[0,1,2]$, & $[0,2,\frac{n+2}{2}]$
\end{tabular}
\end{center}
define a series $A_6(n)$ of cyclic $\{3,6\}$-equivelar triangulations.
The link of the vertex $0$ in $A_{6}(n)$ is one circle consisting of the edges
\begin{center}
\small
\begin{tabular}{c@{\hspace{4mm}}c}
$[1,2]$,     & $[2,\frac{n+2}{2}]$,               \\[1.5mm]
$[1,n-1]$,   & $[\frac{n-2}{2},n-2]$,             \\[1.5mm]
$[n-2,n-1]$, & $[\frac{n-2}{2},\frac{n+2}{2}]$.
\end{tabular}
\end{center}
The Euler characteristic of $A_{6}(n)$ is $\chi=0$, its $f$-vector is $f=(n,3n,2n)$.
For $n=8+4r$, $r\geq 0$, $A_{6}(8+4r)$ is a triangulated torus (see Figure~\ref{fig:A_6_8_4r}),
whereas $A_{6}(10+4r)$ is a triangulated Klein bottle for $r\geq 0$ (see Figure~\ref{fig:A_6_10_4r}).

\subsection{The series \mathversion{bold}$B_{9}(n)$\mathversion{normal}}
\label{subsec:B_9_n}

\begin{figure}
\begin{center}
\small
\psfrag{0}{0}
\psfrag{1}{1}
\psfrag{2}{2}
\psfrag{3}{3}
\psfrag{4}{4}
\psfrag{5}{5}
\psfrag{6}{6}
\psfrag{7}{7}
\psfrag{8}{8}
\psfrag{9}{9}
\psfrag{10}{10}
\psfrag{11}{11}
\psfrag{12}{12}
\psfrag{13}{13}
\psfrag{14}{14}
\psfrag{15}{15}
\includegraphics[width=0.72\linewidth]{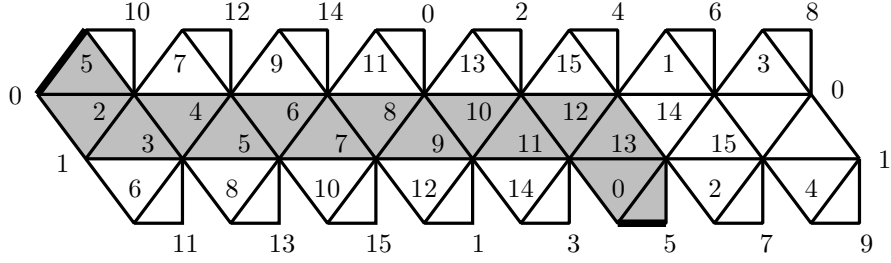}
\end{center}
\caption{The orientable example $B_{9}(16)$.}
\label{fig:B_9_16}
\end{figure}

\begin{figure}
\begin{center}
\small
\psfrag{0}{0}
\psfrag{1}{1}
\psfrag{2}{2}
\psfrag{3}{3}
\psfrag{4}{4}
\psfrag{5}{5}
\psfrag{6}{6}
\psfrag{7}{7}
\psfrag{8}{8}
\psfrag{9}{9}
\psfrag{10}{10}
\psfrag{11}{11}
\psfrag{12}{12}
\psfrag{13}{13}
\psfrag{14}{14}
\psfrag{15}{15}
\psfrag{16}{16}
\psfrag{17}{17}
\psfrag{18}{18}
\psfrag{19}{19}
\includegraphics[width=0.9\linewidth]{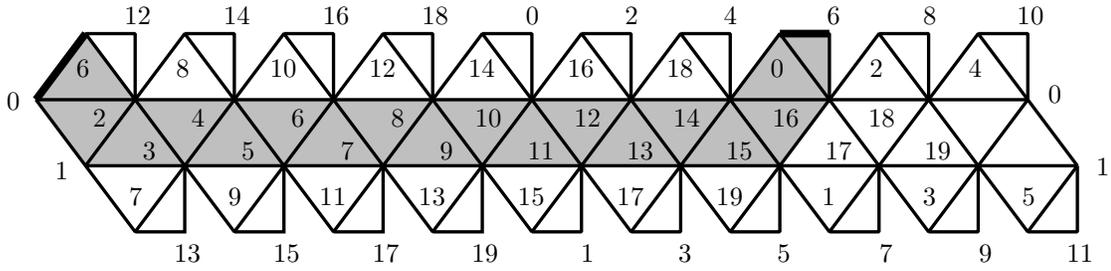}
\end{center}
\caption{The non-orientable example $B_{9}(20)$.}
\label{fig:B_9_20}
\end{figure}

Let $n=16+4r$ for $r\geq 0$. Then the set of generating triangles
\begin{center}
\small
\begin{tabular}{c@{\hspace{4mm}}c@{\hspace{4mm}}c}
$[0,1,2]$, & $[0,2,\frac{n+4}{4}]$, & $[0,\frac{n-4}{4},\frac{n}{2}]$ 
\end{tabular}
\end{center}
defines a series $B_9(n)$ of cyclic $\{3,9\}$-equivelar triangulated surfaces.
The link of vertex $0$ in $B_{9}(n)$ is a circle with nine edges
\begin{center}
\small
\begin{tabular}{c@{\hspace{4mm}}c@{\hspace{4mm}}c}
$[1,2]$,     & $[2,\frac{n+4}{4}]$,               & $[\frac{n-4}{4},\frac{n}{2}]$, \\[1.5mm]
$[1,n-1]$,   & $[\frac{n-4}{4},n-2]$,             & $[\frac{n+4}{4},\frac{3n+4}{4}]$, \\[1.5mm]
$[n-2,n-1]$, & $[\frac{3n-4}{4},\frac{3n+4}{4}]$, & $[\frac{n}{2},\frac{3n-4}{4}]$.
\end{tabular}
\end{center}
For $n=16+8s$, $s\geq 0$, the examples of the series are orientable of genus $g(n)=\frac{n+4}{4}$,
while for $n=20+8s$, $s\geq 0$, the examples are non-orientable of genus $u(n)=\frac{n+4}{2}$.
The Euler characteristic of $B_{9}(n)$ is $\chi=-\frac{n}{2}$, its $f$-vector is $f=(n,\frac{9n}{2},3n)$.
Figure~\ref{fig:B_9_16} displays the orientable example $B_{9}(16)$,
whereas Figure~\ref{fig:B_9_20} displays the non-orientable example $B_{9}(20)$.

\subsection{The series \mathversion{bold}$C_{12}(n)$\mathversion{normal}}
\label{subsec:C_12_n}

For $n=14$, $15$, and $n\geq 17$, the generating triangles 
\begin{center}
\small
\begin{tabular}{c@{\hspace{4mm}}c@{\hspace{4mm}}c@{\hspace{4mm}}c}
$[0,1,2]$, & $[0,2,5]$, & $[0,3,8]$, & $[0,4,8]$
\end{tabular}
\end{center}
define a cyclic series $C_{12}(n)$ of non-orientable triangulated surfaces of genus~$u(n)=n+2$.
The examples of the series are $\{3,12\}$-equivelar, with $C_{12}(n)$ having $f$-vector $f=(n,6n,4n)$
and Euler characteristic $\chi=-n$. 
The link of vertex $0$ in $C_{12}(n)$ is composed of the edges
\begin{center}
\small
\begin{tabular}{c@{\hspace{4mm}}c@{\hspace{4mm}}c@{\hspace{4mm}}c}
$[1,2]$,     & $[2,5]$,     & $[3,8]$,     & $[4,8]$, \\ 
$[1,n-1]$,   & $[3,n-2]$,   & $[5,n-3]$,   & $[4,n-4]$, \\ 
$[n-2,n-1]$, & $[n-5,n-3]$, & $[n-8,n-5]$, & $[n-8,n-4]$.
\end{tabular}
\end{center}

For $n=12$, the above set of generating triangles defines an orientable triangulated surface $C_{10}(12)$ 
with $q=10$ and genus~$g(12)=5$. The first three orbits of this example are of size $12$,
whereas the last orbit is of size $4$. The $f$-vector of $C_{10}(12)$ is $f=(12,60,40)$.

\subsection{The series \mathversion{bold}$D_{12}(n)$\mathversion{normal}}

For $n=13$, $17$, $18$, $19$, and $n\geq 21$, the generating triangles 
\begin{center}
\small
\begin{tabular}{c@{\hspace{4mm}}c@{\hspace{4mm}}c@{\hspace{4mm}}c}
$[0,1,2]$, & $[0,2,6]$, & $[0,4,10]$, & $[0,5,10]$
\end{tabular}
\end{center}
define a cyclic series $D_{12}(n)$ of non-orientable triangulated surfaces of genus~$u(n)=n+2$.
The examples of the series are $\{3,12\}$-equivelar, with $D_{12}(n)$ having $f$-vector $f=(n,6n,4n)$
and Euler characteristic $\chi=-n$. 

For $n=15$, $D_{10}(15)$ is a non-orientable triangulated surface with $q=10$, genus $u(15)=12$
and $f=(15,75,50)$.

\subsection{The series \mathversion{bold}$E_{12}(n)$\mathversion{normal}}
\label{subsec:E_12_n}

Another series $E_{12}(n)$ of non-orientable cyclic $\{3,12\}$-equivelar surfaces
with $f=(n,6n,4n)$ and $\chi=-n$ is defined for $n=17+2m$, $m\geq 0$, 
by the generating triangles
\begin{center}
\small
\begin{tabular}{c@{\hspace{4mm}}c@{\hspace{4mm}}c@{\hspace{4mm}}c}
$[0,1,2]$, & $[0,2,5]$, & $[0,3,\frac{n+3}{2}]$, & $[0,5,\frac{n+5}{2}]$.
\end{tabular}
\end{center}

For $n=15$, this set of generating triangles defines a non-orientable $\{3,10\}$-equivelar surface 
of genus~$u(15)=12$ different from $D_{10}(15)$.
As for $D_{10}(15)$, the first three orbits of $E_{10}(15)$ are of size $15$, 
whereas the last orbit is of size $5$. The $f$-vector of $E_{10}(15)$ is $f=(15,75,50)$.

\subsection{The series \mathversion{bold}$F_{12}(n)$\mathversion{normal}}
\label{subsec:F_12_n}

A further series $F_{12}(n)$ of cyclic $\{3,12\}$-equivelar triangulated surfaces
with $f=(n,6n,4n)$ and $\chi=-n$ is defined for $n=14$ and $n=18+2m$, $m\geq 0$, 
by the generating triangles
\begin{center}
\small
\begin{tabular}{c@{\hspace{4mm}}c@{\hspace{4mm}}c@{\hspace{4mm}}c}
$[0,1,2]$, & $[0,2,6]$, & $[0,3,6]$, & $[0,4,\frac{n+4}{2}]$.
\end{tabular}
\end{center}

For $n=14+4r$, $r\geq 0$, the examples $F_{12}(n)$ are orientable
of genus $g(n)=\frac{n+2}{2}$, while for $n=20+4r$, $r\geq 0$, the examples $F_{12}(n)$ 
are non-orientable of genus $u(n)=n+2$.

\subsection{The series \mathversion{bold}$G_{18}(n)$ and $H_{18}(n)$\mathversion{normal}}
\label{subsec:G_18_n_H_18_n}

For $n=19+2m$, $m\geq 0$, there is an orientable cyclic series of $\{3,18\}$-equivelar 
triangulated surfaces~$G_{18}(n)$ and a non-orientable cyclic series of $\{3,18\}$-equivelar 
triangulated surfaces $H_{18}(n)$, both with $f=(n,9n,6n)$ and $\chi=-2n$.

\begin{figure}
\begin{center}
\scriptsize
\psfrag{0}{0}
\psfrag{1}{1}
\psfrag{2}{2}
\psfrag{3}{3}
\psfrag{4}{4}
\psfrag{5}{5}
\psfrag{6}{6}
\psfrag{7}{7}
\psfrag{8}{8}
\psfrag{9}{9}
\psfrag{10}{10}
\psfrag{11}{11}
\psfrag{12}{12}
\psfrag{13}{13}
\psfrag{14}{14}
\psfrag{15}{15}
\psfrag{16}{16}
\psfrag{17}{17}
\psfrag{18}{18}
\includegraphics[width=\linewidth]{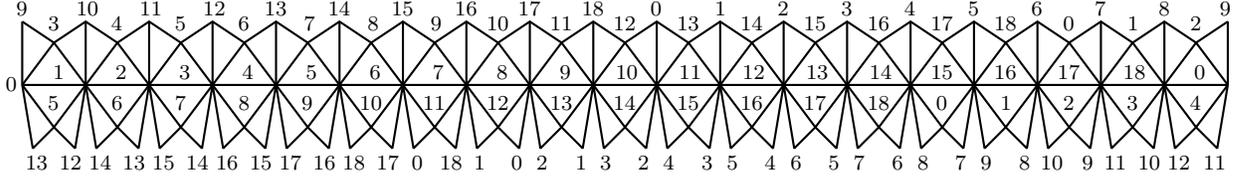}
\end{center}
\caption{The orientable example $G_{18}(19)$.}
\label{fig:G_18_19}
\end{figure}

The examples of the series $G_{18}(n)$ have genus $g(n)=n+1$ and are defined by the generating triangles
\begin{center}
\small
\begin{tabular}{c@{\hspace{4mm}}c@{\hspace{4mm}}c@{\hspace{4mm}}c@{\hspace{4mm}}c@{\hspace{4mm}}c}
$[0,1,3]$, & $[0,1,5]$, & $[0,2,\frac{n-1}{2}]$, & $[0,3,\frac{n-1}{2}]$, & $[0,4,\frac{n+3}{2}]$, & $[0,5,\frac{n+7}{2}]$.
\end{tabular}
\end{center}
The example $G_{18}(19)$ is displayed in Figure~\ref{fig:G_18_19}.

The examples of the series $H_{18}(n)$ have genus $u(n)=2n+2$ and are defined by the generating triangles
\begin{center}
\small
\begin{tabular}{c@{\hspace{4mm}}c@{\hspace{4mm}}c@{\hspace{4mm}}c@{\hspace{4mm}}c@{\hspace{4mm}}c}
$[0,1,2]$, & $[0,2,5]$, & $[0,3,\frac{n+1}{2}]$, & $[0,4,\frac{n+1}{2}]$, & $[0,4,\frac{n+3}{2}]$, & $[0,5,\frac{n+7}{2}]$.
\end{tabular}
\end{center}

\subsection{The series \mathversion{bold}$I_{18}(n)$\mathversion{normal}}

The cyclic $\{3,18\}$-equivelar series $I_{18}(n)$ with $f=(n,9n,6n)$ and $\chi=-2n$ is defined 
for $n=20$ and $n=24+2m$, $m\geq 0$, by the set of generating triangles
\begin{center}
\small
\begin{tabular}{c@{\hspace{4mm}}c@{\hspace{4mm}}c@{\hspace{4mm}}c@{\hspace{4mm}}c@{\hspace{4mm}}c}
$[0,1,2]$, & $[0,2,7]$, & $[0,3,7]$, & $[0,3,\frac{n-2}{2}]$, & $[0,4,\frac{n+4}{2}]$, & $[0,5,\frac{n+8}{2}]$.
\end{tabular}
\end{center}

For $n=20+4r$, $r\geq 0$, the examples $I_{18}(n)$ 
are non-orientable of genus $u(n)=2n+2$, while the examples $I_{18}(n)$ 
are orientable of genus $g(n)=n+1$ for $n=26+4r$, $r\geq 0$.

\subsection{Orientability versus non-orientability}
\label{subsec:or_vs_nonor}

The series  $A_{6}(n)$, $B_{9}(n)$, $C_{12}(n)$, $D_{12}(n)$, $E_{12}(n)$, $F_{12}(n)$, $H_{18}(n)$, 
and $I_{18}(n)$ all contain the orbit of the generating triangle $[0,1,2]$. 
This orbit is a triangulated cylinder if $n$ is even and is a triangulated M\"obius strip 
whenever $n$ is odd. It therefore immediately follows that the examples of the series $E_{12}(n)$ and $H_{18}(n)$
and all the examples of the subseries $C_{12}(15+2m)$, $m\geq 0$,
and of the subseries $D_{12}(13+2m)$, $m\geq 0$, contain a M\"obius strip and hence are non-orientable.

We next inspect the triangulations $B_{9}(16)$ and $B_{9}(20)$, which each are composed of
three orbits of $16$ respectively $20$ triangles. In the Figures~\ref{fig:B_9_16} and \ref{fig:B_9_20} 
the cylinder defined by the orbit of the triangle $[0,1,2]$ lies horizontally. 
Attached to it are the triangles of the orbit $[0,2,\frac{n+4}{4}]$ and finally the
triangles of the orbit $[0,\frac{n-4}{4},\frac{n}{2}]$, for $n=16$ and $n=20$, respectively.
In both cases (Figures~\ref{fig:B_9_16} and \ref{fig:B_9_20}), we obtain a triangulated disk with identified edges on the boundary.
We start at vertex $0$ at the left and follow the boundary of the triangulated disk 
in clockwise direction. In $B_{9}(20)$, the edge $[0,6]$ appears twice with the same orientation, 
which forces a M\"obius strip (in grey). Hence, $B_{9}(20)$ is non-orientable. In contrast, 
identified boundary edges in $B_{9}(16)$ appear with opposite orientations (so that, in particular, 
the respective grey band in between the two copies of the edge $[0,5]$ is a cylinder). It follows 
that $B_{9}(16)$ is orientable. 

It is straight forward to generalize the orientability discussion of the cases $B_{9}(16)$ and $B_{9}(20)$
to conclude that all examples of the subseries  $B_{9}(16+8s)$, $s\geq 0$, are orientable,
whereas the examples of the subseries  $B_{9}(20+8s)$, $s\geq 0$, are non-orientable.
Likewise, respective orientability results follow for $A_{6}(n)$, $C_{12}(n)$, $D_{12}(n)$, $F_{12}(n)$, 
and $I_{18}(n)$.

Finally, Figure~\ref{fig:G_18_19} shows that the triangulation $G_{18}(19)$ is orientable, 
as are all the examples of the series $G_{18}(19+2m)$, $m\geq 0$.

\section{The \mathversion{bold}$[0,1,2]$-family of cyclic triangulations\mathversion{normal}}
\label{sec:012}

We are now going to describe a large class of cyclic triangulations of (closed) surfaces. 
In particular, we will prove (in Section~\ref{sec:012_series}) our Main Theorem 
that for every $q=3k$, $k\geq 2$, and every $q=3k+1$, $k\geq 3$, there are infinitely 
many orientable and non-orientable $q$-equivelar cyclic triangulations of surfaces.

\begin{deff}
Let the cyclic group action on the $n$ vertices $0,1,\dots,n-1$ be defined by the cyclic shift $(0,1,\dots,n-1)$.
The \emph{$[0,1,2]$-family of cyclic triangulations} is the family of all cyclic triangulations of surfaces
that contain the orbit of the triangle $[0,1,2]$.
\end{deff}

We will see in the following, that the fundamental domains of the members of the $[0,1,2]$-family
have a particularly nice description that allows to easily read off the orientability or non-orientability
of various subseries.

As we noted before, the orbit of the triangle $[0,1,2]$ is a cylinder if $n$ is even and is a M\"obius band 
whenever $n$ is odd. The edge $[0,1]$ and all its translates are interior edges of the $[0,1,2]$-orbit, 
whereas the edge $[0,2]$ and its translates are boundary edges. If we cut open the $[0,1,2]$-band along 
the edge $[0,1]$, we obtain a triangulated disk with identifications on the boundary: 
the boundary of the disk contains the two copies of the edge $[0,1]$, whereas all other edges 
of the boundary are the cyclic translates of the edge $[0,2]$. 
By the strong connectivity of the triangulated surfaces of the $[0,1,2]$-family, 
the edge $[0,2]$ lies in a second triangle $[0,2,x]$. As in the Figures~\ref{fig:A_6_8_4r}--\ref{fig:B_9_20},
we glue the triangle $[0,2,x]$ and its translates to the boundary edges of the cut open $[0,1,2]$-band.
The resulting triangulated disk may have all its boundary edges pairwise identified; 
see for example Figures~\ref{fig:A_6_8_4r}--\ref{fig:A_6_10_4r}.
In case the edge $[0,x]$ (respectively the edge $[2,x]$) is, at this stage, not identified with another boundary edge, 
then $[0,x]$ (respectively the edge $[2,x]$) is contained in a triangle $[0,x,y]$ (respectively $[2,x,y]$),
which generates a third orbit of triangles. The triangles of this orbit are again glued to the triangulated disk.
We proceed iteratively and add further triangles and translates thereof until eventually all boundary edges of 
the obtained triangulated disk are pairwise identified. In the resulting triangulated disk
a tree of orbit generating triangles is glued to the triangle $[0,1,2]$ along the edge $[0,2]$.

\begin{deff}
Let a \emph{seed} of a member of the $[0,1,2]$-family be a fundamental domain consisting
of the triangle $[0,1,2]$ and a tree of distinct orbit generating triangles that are glued 
to the triangle $[0,1,2]$ of the cut open $[0,1,2]$-band along the edge~$[0,2]$.
\end{deff}

For the members of the series $A_6(n)$ (Figures~\ref{fig:A_6_8_4r}--\ref{fig:A_6_10_4r}), 
the seed consists of the two triangles $[0,1,2]$ and $[0,2,\frac{n+2}{2}]$.
For the examples of the series $B_9(n)$, a seed is formed by the three triangles $[0,1,2]$, $[0,2,\frac{n+4}{4}]$, 
and $[2,\frac{n+4}{4},\frac{n+4}{2}]$. There are two ways to glue the triangle $[2,\frac{n+4}{4},\frac{n+4}{2}]$
to the previous two triangles (along the edge $[0,\frac{n+4}{4}]$ or along the edge $[2,\frac{n+4}{4}]$).
Thus, in general, a member of the $[0,1,2]$-family can have different \emph{seed trees}.

\subsubsection*{The standard case:} 

If for a member of the $[0,1,2]$-family all orbits of triangles are of size $n$,
then a seed has the following properties:
\begin{itemize}
\itemsep=0pt
\item[S0.] Each triangle of a seed has three distinct vertices and appears only once in the seed.
           The triangle $[0,1,2]$ is a leaf of every corresponding seed tree.
\item[S1.] If a seed (or more precise, a seed tree of triangles) has an even number of boundary edges, 
  then its boundary edges can be grouped pairwise, as described in S3.
\item[S2.] If a seed has an odd number of boundary edges, then $n$ is even and exactly 
  one of the boundary edges is a translate of the edge $[0,\frac{n}{2}]$, whereas
  all other of its boundary edges can be grouped pairwise, as described in S3.
\item[S3.] Two boundary edges of a seed form a pair if they 
  are cyclic images of each other. If a boundary edge $[a,b]$ is not a
  translate of the edge $[0,\frac{n}{2}]$, then to $[a,b]$ there is exactly 
  one distinct cyclic image on the boundary of the seed; no cyclic 
  image of $[a,b]$ is an interior edge of the seed.
\item[S4.] An interior edge is not the cyclic image of another interior
  edge or of a boundary edge.
\item[S5.] Since seeds have a tree structure, a seed has an odd (even) number
  of boundary edges exactly when it has an odd (even) number of triangles.
\end{itemize}

If we glue together all the cyclic images of a seed along the edges $[1,2]$, $[2,3]$, \dots, $[0,n-1]$, 
we obtain a triangulated disk with identifications on the boundary. 
\begin{lem}
If a tree of orbit generating triangles has properties S0--S5, then the disk 
resulting from glueing all cyclic images of the tree along the edges $[1,2]$, $[2,3]$, \dots, $[0,n-1]$ 
defines a cyclic triangulation of a surface or a pinched surface
(a surface with isolated singularities).
\end{lem}

\textbf{Proof:} By S0, the collection of triangles obtained by uniting the cyclic images
of the tree of orbit generating triangles is a simplicial complex $K$ with cyclic symmetry.
According to \mbox{S1--S4}, every edge of $K$ is contained in exactly two triangles, 
which implies that $K$ is a triangulated surface or a triangulated pinched surface.~\hfill$\Box$

\subsubsection*{The exceptional case:}
 
It remains to consider the case that not all orbits of triangles of a member of the $[0,1,2]$-family 
are of size $n$. This can only happen when $n$ is divisible by $3$, with the only cyclic orbit of triangles 
of size less than $n$ being the orbit generated by the triangle $[0,\frac{n}{3},\frac{2n}{3}]$.
This orbit has $\frac{n}{3}$ triangles and is bounded by the $n$ edges of the orbit of the edge $[0,\frac{n}{3}]$. 
In particular, there is a unique orbit of triangles neighboring to the orbit generated by the triangle $[0,\frac{n}{3},\frac{2n}{3}]$.
For the seed, it follows that the respective translate of the triangle $[0,\frac{n}{3},\frac{2n}{3}]$
used in the seed is a leaf of the seed tree. 

\begin{deff}
Let $n$ be divisible by $3$. If a member of the $[0,1,2]$-family contains the orbit generated 
by the triangle $[0,\frac{n}{3},\frac{2n}{3}]$, then the corresponding \emph{reduced seed}
is the tree consisting of the triangle $[0,1,2]$ and the tree of distinct orbit generating
triangles of size $n$ that are glued to the triangle $[0,1,2]$ of the cut open $[0,1,2]$-band 
along the edge~$[0,2]$.
\end{deff}

Any reduced seed has the property
\begin{itemize}
\itemsep=0pt
\item[S$\frac{1}{3}$.] The boundary of a reduced seed contains exactly one translate of the edge $[0,\frac{n}{3}]$
  and is otherwise described by the rules S0--S5, where in S1--S5 the numbers of boundary edges have to be increased by one,
  respectively.
\end{itemize}

The observation that every member of the $[0,1,2]$-family contains at most one orbit of size less than~$n$
is a property that holds for all cyclic triangulations of surfaces.

\begin{prop}
A cyclic $n$-vertex triangulation of a surface contains at most one orbit of size less than $n$.
If such a short orbit is present, then $n$ is divisible by $3$ and the short orbit is generated 
by the triangle $[0,\frac{n}{3},\frac{2n}{3}]$.
\end{prop}

\begin{cor}\label{cor:cyclic}
Cyclic triangulations of surfaces either are $q$-equivelar with $q=3k$ for some $k\geq 1$
or are $q$-equivelar with $q=3k+1$ for some $k\geq 1$.
\end{cor}

The requirements imposed by the rules S0--S5, respectively S$\frac{1}{3}$ in the presence of a short orbit,
are moderate and allow for an abundance of constructions that lead to infinite series 
of cyclic triangulations of surfaces (and pinched surfaces). In the next section, 
we present various such series of examples and discuss the $q=7$. Hereby, we prove our Main Theorem.

\section{Some series of cyclic \mathversion{bold}$[0,1,2]$-triangulations\mathversion{normal}}
\label{sec:012_series}

\subsection{The series \mathversion{bold}$S_{k,n}$\mathversion{normal}}
\label{subsec:S_k_n}

\begin{figure}
\begin{center}
\small
\psfrag{0}{0}
\psfrag{1}{1}
\psfrag{2}{2}
\psfrag{5}{5}
\psfrag{10}{10}
\psfrag{16}{16}
\psfrag{22}{22}
\psfrag{29}{29}
\psfrag{36}{36}
\psfrag{19}{19}
\psfrag{[1]}{\scriptsize $\langle 1\rangle$}
\psfrag{[2]}{\scriptsize $\langle 2\rangle$}
\psfrag{[3]}{\scriptsize $\langle 3\rangle$}
\psfrag{[5]}{\scriptsize $\langle 5\rangle$}
\psfrag{[6]}{\scriptsize $\langle 6\rangle$}
\psfrag{[7]}{\scriptsize $\langle 7\rangle$}
\psfrag{[17]}{\scriptsize $\langle 17\rangle$}
\psfrag{[8]}{\scriptsize $\langle 8\rangle$}
\psfrag{[14]}{\scriptsize $\langle 14\rangle$}
\psfrag{[20]}{\scriptsize $\langle 20\rangle$}
\psfrag{[27]}{\scriptsize $\langle 27\rangle$}
\psfrag{[34]}{\scriptsize $\langle 34\rangle$}
\includegraphics[width=0.675\linewidth]{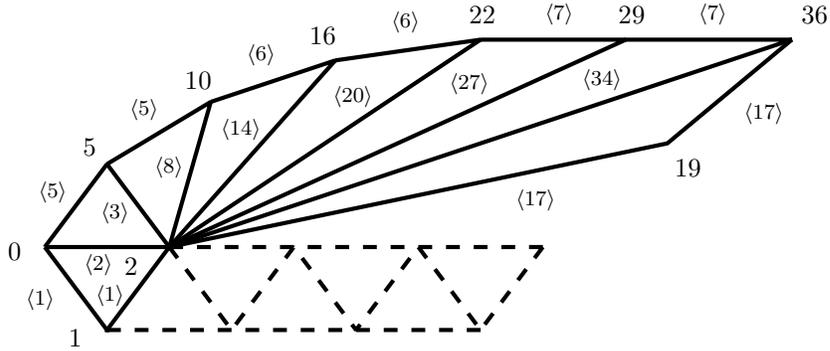}
\end{center}
\caption{The seed of $S_{4,n}$.}
\label{fig:S_4_n}
\end{figure}

We define a two-parameter family of cyclic triangulations $S_{k,n}$: For every $k\geq 2$, let 
the seed of~$S_{k,n}$ be composed of the $2k$ triangles

\medskip

\begin{tabular}{l}
$[0,1,2]$, \\[2mm]
$[0,2,5]$,\, $[2,5,10]$, \\[2mm]
$[2,10,16]$,\, $[2,16,22]$, \\[2mm]
\dots, \\[2mm]
$[2,2\displaystyle\sum_{i=2}^{k-1}(i+3),2\displaystyle\sum_{i=2}^k(i+3)-(k+3)]$,\, $[2,2\displaystyle\sum_{i=2}^k(i+3)-(k+3),2\displaystyle\sum_{i=2}^k(i+3)]$, \\
$[2,2+\frac{1}{2}(2\displaystyle\sum_{i=2}^k(i+3)-2),2\displaystyle\sum_{i=2}^k(i+3)]$. 
\end{tabular}

\medskip

The subseries $S_{2,n}$ of $S_{k,n}$ then coincides with the one-parameter series $C_{12}(n)$. 
The seed of the cyclic triangulations~$S_{4,n}$ is displayed in Figure~\ref{fig:S_4_n}.

Since all the triangles of the seeds of $S_{2,n}$, $S_{3,n}$, \dots  contain the vertex $2$, 
a more compressed way to represent the seed of $S_{k,n}$ is by giving the ``lower'' boundary edges
$[0,1]$ and $[1,2]$ as well as the ``upper'' boundary edges, noted as a path
$$0\stackrel{\langle 5\rangle}{\mbox{---}}5\stackrel{\langle 5\rangle}{\mbox{---}}10
\stackrel{\langle 6\rangle}{\mbox{---}}16\stackrel{\langle 6\rangle}{\mbox{---}}22\mbox{---}\,\dots\,\mbox{---}2\displaystyle\sum_{i=2}^{k-1}(i+3)
\stackrel{\langle k+3\rangle}{\mbox{---}}\left(2\displaystyle\sum_{i=2}^k(i+3)-(k+3)\right)\stackrel{\langle k+3\rangle}{\mbox{---}}2\displaystyle\sum_{i=2}^k(i+3)$$
$$\stackrel{\langle -\frac{1}{2}(2\sum_{i=2}^k(i+3)-2)\rangle}{\mbox{---}}2+\frac{1}{2}(2\displaystyle\sum_{i=2}^k(i+3)-2)\stackrel{\langle -\frac{1}{2}(2\sum_{i=2}^k(i+3)-2)\rangle}{\mbox{---}}2$$
or simply by giving the differences $\langle 5\rangle$, $\langle 5\rangle$, $\langle 6\rangle$, $\langle 6\rangle$, \dots, $\langle k+3\rangle$, $\langle k+3\rangle$,
$\langle-\frac{1}{2}(2\sum_{i=2}^k(i+3)-2)\rangle$, $\langle-\frac{1}{2}(2\sum_{i=2}^k(i+3)-2)\rangle$. 
These differences are added successively to the initial vertex $0$ until we reach the vertex $2$.
The first $2(k-1)$ differences define $2(k-1)$ orbit generating triangles of a respective seed, 
whereas the last two differences $\langle\frac{1}{2}(2\sum_{i=2}^k(i+3)-2)\rangle=\langle\frac{1}{2}(k^2+7k)-5\rangle$
represent two boundary edges of the orbit generating triangle $[2,2+\frac{1}{2}(2\sum_{i=2}^k(i+3)-2),2\sum_{i=2}^k(i+3)]$.

For $k=1$, we can define the seed of $S_{1,n}$ by gluing the triangle $[2,2+\frac{1}{2}(n-2),n]$
to the triangle $[0,1,2]$ along the edge $[2,n]=[0,2]$. The series $S_{1,n}$ then coincides with
the series~$A_{6}(n)$.

\subsubsection*{The link \mathversion{bold}of $0$ in $S_{k,n}$ for $k\geq 2$:\mathversion{normal}}

\begin{figure}
\begin{center}
\small
\psfrag{0}{0}
\psfrag{1}{1}
\psfrag{-1}{$-1$}
\psfrag{2}{2}
\psfrag{-2}{$-2$}
\psfrag{3}{3}
\psfrag{-3}{$-3$}
\psfrag{5}{5}
\psfrag{-5}{$-5$}
\psfrag{6}{6}
\psfrag{-6}{$-6$}
\psfrag{8}{8}
\psfrag{-8}{$-8$}
\psfrag{10}{10}
\psfrag{-10}{$-10$}
\psfrag{14}{14}
\psfrag{-14}{$-14$}
\psfrag{20}{20}
\psfrag{-20}{$-20$}
\includegraphics[width=0.25\linewidth]{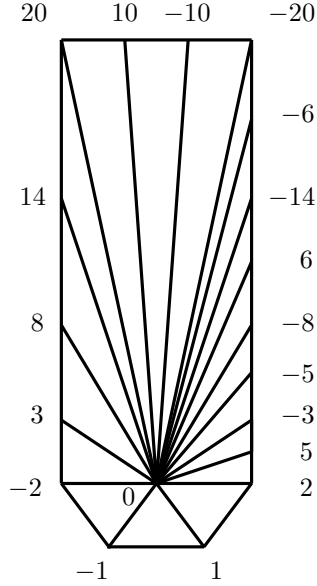}
\end{center}
\caption{The link of $0$ in $S_{3,n}$.}
\label{fig:S_3_n_link}
\end{figure}

The link of $0$ in $S_{k,n}$ can easily be read off the seed of $S_{k,n}$:
For every vertex $v$ of the seed, take all triangles of the seed that contain $v$
and translate these triangles by $-v$. The union over all these triangles is the link of $0$ in $S_{k,n}$.
It immediately follows that for every distance $\langle d\rangle$ in the seed, $d$ and $-d$ are vertices
of the link of $0$. Also, for a neighbor $w$ of $2$ in the seed, $w-2$ and $-w+2$ are vertices of the link of $0$.
Figure~\ref{fig:S_3_n_link} displays the link of $0$ in $S_{3,n}$. This symmetric depiction of the link of $0$ 
in $S_{3,n}$ can be extended to a description of the link of $0$ in $S_{k,n}$ for all $k\geq 2$. In particular, 
if there are no forbidden identifications (see below) in the case of small $n$, then the link of $0$ in $S_{k,n}$ is 
a triangulated disk and thus $S_{k,n}$ is a triangulated surface.

\subsubsection*{Range of \mathversion{bold}$n$:\mathversion{normal}} 

According to Section~\ref{subsec:C_12_n}, the series $C_{12}(n)=S_{2,n}$
is defined for $n=14$, $15$, and $n\geq 17$. For $k=4$, the series $S_{4,n}$ yields triangulated surfaces 
for $n=38$, $43$, $45$, $46$, $49$, $50$, $52$, $53$, $55$, \dots,~$60$, $62$, \dots,~$67$,
and $n\geq 69$. In these cases, the respective examples have $f$-vector $f=(n,12n,8n)$ and Euler characteristic $\chi=-3n$.
In addition, $S_{4,51}$ is a triangulated surface with a short orbit generated by the triangle $[0,17,34]$; 
its $f$-vector is $f=(51,561,374)$, its Euler characteristic $\chi=-136$. For $n=12$, $C_{10}(12)=S_{2,12}$
has a short orbit of size $4$.

If $n=68$, then the orbit of the edge $2\stackrel{\langle 34\rangle}{\mbox{---}}36$
is a short orbit of size $34$, which causes forbidden identifications of interior edges.
Similarly, if $n=54$, then the orbit of the edge $2\stackrel{\langle 27\rangle}{\mbox{---}}29$
is a short forbidden orbit of size $27$.

If $n=61$, then the orbit generated by the edge $2\stackrel{\langle 27\rangle}{\mbox{---}}29$
coincides with the orbit generated by the edge $2\stackrel{\langle 34\rangle}{\mbox{---}}36$,
since $27+34=61$; thus $n=61$ is not permitted.

We conclude that $S_{k,n}$ \emph{is not} a triangulated surface if $n$ is the double of some occurring difference 
or the sum of two distinct occurring differences. If we let $d=2\sum_{i=2}^k(i+3)-2$, 
then $S_{k,n}$ is a triangulated surface for all $n\geq 2d+1$. The smallest admissible value for $n$ is $n=d+4$.
Moreover, the value $n=\frac{3d}{2}$ is admissible, in which case $S_{k,3\sum_{i=2}^k(i+3)-3}$ has a small orbit.
 
\subsubsection*{Non-orientability of \mathversion{bold}$S_{k,n}$ for $k\geq 2$:\mathversion{normal}} 

If $n$ is odd, then, as remarked in Section~\ref{subsec:or_vs_nonor}, $S_{k,n}$ is non-orientable,
since then the orbit of the triangle $[0,1,2]$ is a M\"obius band in $S_{k,n}$.
If $n$ is even, then the orbit of $[0,1,2]$ is a cylinder. However, as we have seen in Figure~\ref{fig:B_9_20}, 
other M\"obius bands might be present.

In fact, for every boundary edge of the seed it has to be examined whether the edge and its second copy on
the boundary of the triangulated disk generated by the seed yields a M\"obius band (as in Figure~\ref{fig:B_9_20}
for the edge $[0,6]$) or a cylinder (as in Figure~\ref{fig:B_9_16} for the edge $[0,5]$).
For even $n$, the cylinder generated by the triangle $[0,1,2]$ splits the triangulated disk
into an upper and a lower part. Let both copies of a boundary edge have the same orientation. 
If the two copies both lie in the upper part, they span a M\"obius band. If the second copy lies in the lower part, 
a cylinder is induced. The triangles of the upper part are glued to the central cylinder via the ``even'' edges
$[0,2]$, $[2,4]$, \dots, $[n-2,0]$, while the triangles of the lower part are glued to the central 
cylinder via the ``odd'' edges $[1,3]$, $[3,5]$, \dots, $[n-1,1]$. It follows, for example, that for 
the pair of edges $0\stackrel{\langle 5\rangle}{\mbox{---}}5$ and $5\stackrel{\langle 5\rangle}{\mbox{---}}10$
of difference $5$ in the seed of $S_{4,n}$ the second copy of $0\stackrel{\langle 5\rangle}{\mbox{---}}5$, 
as a translate of $5\stackrel{\langle 5\rangle}{\mbox{---}}10$ by the odd number $-5$, lies in the lower part of the triangulated disk, 
thus inducing a cylinder. 
In contrast, for the pair of edges $10\stackrel{\langle 6\rangle}{\mbox{---}}16$ and $16\stackrel{\langle 6\rangle}{\mbox{---}}22$
of difference $6$ in the seed of $S_{4,n}$ the second copy of $10\stackrel{\langle 6\rangle}{\mbox{---}}16$, 
as a translate of $16\stackrel{\langle 6\rangle}{\mbox{---}}22$ by $-6$, lies in the upper part of the triangulated disk, 
thus inducing a M\"obius band. 
We conclude that $S_{k,n}$ is non-orientable for all admissible $n$ if $k\geq 3$.
For $k=2$, the edge $10\stackrel{\langle -4\rangle}{\mbox{---}}6$ and its copy, as a translate
of the edge $6\stackrel{\langle -4\rangle}{\mbox{---}}2$ by $4$, both lie in the upper part,
which implies the non-orientability of the examples $S_{2,n}$ for all admissible~$n$.

For $n=8+4r$, $r\geq 0$, the examples $S_{1,n}$ are triangulated tori (and thus orientable), whereas 
$S_{1,10+4r}$ is a triangulated Klein bottle for $r\geq 0$; see Section~\ref{subsec:A_6_n}.

\begin{thm}
Let $q=6k$ for $k\geq 1$ and let $n\geq 2(2\sum_{i=2}^k(i+3)-2)+1$. Then the series $S_{k,n}$
provides infinitely many non-orientable cyclic triangulations, which are equivelar of type $\{3,q\}=\{3,6k\}$.
\end{thm}

\subsection{The families \mathversion{bold}$G_{d_1,d_1,d_2,d_2,\dots,d_{k-1},d_{k-1}}$\mathversion{normal} and \mathversion{bold}$\tilde{G}_{d_1,d_1,d_2,d_2,\dots,d_{k-1},d_{k-1}}$\mathversion{normal}}
\label{sec:G}

Instead of the sequence of differences 
$\langle 5\rangle$, $\langle 5\rangle$, $\langle 6\rangle$, $\langle 6\rangle$, \dots, $\langle k+3\rangle$, $\langle k+3\rangle$,
which was used for $S_{k,n}$ above, we can easily modify the construction of the series by admitting arbitrary 
distinct differences $d_1, \dots, d_{k-1}\geq 5$.
For given $k$, any member of the family $G_{d_1,d_1,d_2,d_2,\dots,d_{k-1},d_{k-1}}$, defined by using 
each of the differences $d_1, \dots, d_{k-1}$ and the closing difference $-\frac{1}{2}(2\sum_{i=1}^{k-1}d_i-2)$ 
exactly twice, is a triangulated (pinched) surface for all $n\geq 2(2\sum_{i=1}^{k-1}d_i-2)+1$
(as well as for some smaller values of $n$). If care is taken to avoid non-admissible identifications
of edges, then also negative differences as well as differences $\langle 3\rangle$ and $\langle 4\rangle$ 
can be used to define seeds for respective series of cyclic triangulations; see below and Section~\ref{sec:tessellations} for explicit examples. 
Let the family $\tilde{G}_{d_1,d_1,d_2,d_2,\dots,d_{k-1},d_{k-1}}$ be the collection of such general examples.
All the members of the family $\tilde{G}_{d_1,d_1,d_2,d_2,\dots,d_{k-1},d_{k-1}}$ 
have a pathlike seed of triangles.

\subsection{The pinched series \mathversion{bold}$P_{k,n}$\mathversion{normal}}

Let for odd $k\geq 3$ the examples $P_{k,n}$ be defined by the seed given by the differences 
$\langle 5\rangle$, $\langle 6\rangle$, $\langle 5\rangle$, $\langle 6\rangle$, $\langle 7\rangle$, $\langle 8\rangle$, $\langle 7\rangle$, $\langle 8\rangle$, 
\dots, $\langle k+2\rangle$, $\langle k+3\rangle$, $\langle k+2\rangle$, $\langle k+3\rangle$
and completed by the triangle $[2,2+\frac{1}{2}(2\displaystyle\sum_{i=2}^k(i+3)-2),2\displaystyle\sum_{i=2}^k(i+3)]=[2,2+\frac{1}{2}(k^2+7k)-5,k^2+7k-8]$.

By alternating odd and even differences, $\langle i\rangle$, $\langle i+1\rangle$, $\langle i\rangle$, $\langle i+1\rangle$, 
each pair has one edge in the upper part and the corresponding copy in the lower part. 
Let $k=2\ell+1$ for $\ell\geq 1$, then $\frac{1}{2}(k^2+7k)-5=2\ell^2+9\ell-1$,
which is even for odd $\ell$ and odd for even $\ell$. 
Thus, for even $\ell\geq 2$, the edges $[2,2+\frac{1}{2}(k^2+7k)-5]$ and $[2+\frac{1}{2}(k^2+7k)-5,k^2+7k-8]$
of the closing triangle have their copies in the lower part, whereas for odd $\ell\geq 1$, 
the respective copies are in the upper part.

The examples $P_{k,n}$ are pinched triangulations: If we shift the two pairs of triangles $[0,2,5]$, $[2,5,11]$ and $[2,11,16]$, $[2,16,22]$ of the seed 
by $-5$ and $-16$, respectively, then we obtain a cone $[-5,-3,0]$, $[-3,0,6]$, $[-14,-5,0]$, $[-14,0,6]$ with apex $0$ and its link
the circle $[-5,-3]$, $[-3,6]$, $[6,-14]$ $[-14,-5]$. Thus, the link of $0$ has more than one component.

\begin{thm}
Let $q=6k$ for odd $k\geq 1$ and let $n\geq 2(2\sum_{i=2}^k(i+3)-2)+1$. Then the series $P_{k,n}$
provides infinitely many cyclic triangulations of pinched surfaces, which are equivelar of type $\{3,q\}=\{3,6k\}$.
For $k=2\ell+1$, $\ell\geq 1$, the examples $P_{2\ell+1,n}$ are orientable if $\ell$ is even and non-orientable
if $\ell$ is odd.
\end{thm}

\subsection{The series \mathversion{bold}$T_{k,n}$\mathversion{normal}}

We next suitably modify the seed of $S_{k,n}$ for odd $k=2\ell+1$, $\ell\geq 1$, to obtain a series of orientable cyclic triangulations
of closed surfaces, which are equivelar of type $\{3,q\}=\{3,6+12\ell\}$.
For $k=2\ell+1$, $\ell\geq 1$, we consider the subfamily $T_{k,n}$ of $G_{d_1,d_1,d_2,d_2,\dots,d_{k-1},d_{k-1}}$,
defined by the odd differences $\langle 5\rangle$, $\langle 5\rangle$, $\langle 7\rangle$, $\langle 7\rangle$, 
\dots, $\langle 2k-1\rangle$, $\langle 2k-1\rangle$, $\langle 2k+1\rangle$, $\langle 2k+1\rangle$.
By using odd differences, these pairs have one edge in the upper part and the corresponding copy in the lower part. 
The seed is completed by the triangle $[2,2+\frac{1}{2}(2\sum_{i=2}^k(2i+1)-2),2\sum_{i=2}^k(2i+1)]=[2,k^2+2k-2,2k^2+4k-6]$.
Since $\frac{1}{2}(2\sum_{i=2}^k(2i+1)-2)=k^2+2k-4$ is odd for odd $k$, the respective copies of the two closing edges
lie in the lower part

For $k=1$, let $T_{1,n}=S_{1,n}=A_6(n)$. Then for $n=8+4r$, $r\geq 0$, the examples $T_{1,n}$ are triangulated tori.

\begin{thm}
Let $q=6+12\ell$ for $k=2\ell+1$, $\ell\geq 1$, and let $n\geq 2(2\sum_{i=2}^k(2i+1)-2)+2$ be even. Then the series $T_{k,n}$
provides infinitely many orientable cyclic triangulations, which are equivelar of type $\{3,q\}=\{3,6+12\ell\}$, $\ell\geq 1$.
For $q=6$, the series $T_{1,8+4r}$, $r\geq 0$, provides infinitely many triangulated tori which are equivelar of type $\{3,6\}$.
\end{thm}
 
The generalized Ringel series from \cite{LutzSulankeTiwariUpadhyay2010pre} provides for all $n\geq 7+12\ell$ yet another 
two-parameter series of cyclic orientable triangulations with $q=6+12\ell$, $\ell\geq 0$.

\enlargethispage*{10mm}

\subsection{The series \mathversion{bold}$\overline{T}_{k,n}$\mathversion{normal}}

Another modification of $S_{k,n}$ leads to a series $\overline{T}_{k,n}$ of orientable cyclic triangulations, 
which are equivelar of type $\{3,q\}=\{3,6k\}$ for $k\geq 1$.
As before for $T_{k,n}$, we again use the odd differences $\langle 5\rangle$, $\langle 5\rangle$, $\langle 7\rangle$, $\langle 7\rangle$, 
\dots, $\langle 2k-1\rangle$, $\langle 2k-1\rangle$, $\langle 2k+1\rangle$, $\langle 2k+1\rangle$,
but this time we complete the seed by the triangle 
$[2,2+\frac{1}{2}(2\sum_{i=2}^k(2i+1)-2)+\frac{n}{2},2\sum_{i=2}^k(2i+1)]=[2,k^2+2k-2+\frac{n}{2},2k^2+4k-6]$.
For odd $k$, the difference $\frac{1}{2}(2\sum_{i=2}^k(2i+1)-2)+\frac{n}{2}=k^2+2k-4+\frac{n}{2}$ is odd for $n=4r$,
so that the respective copies of the two closing edges lie in the lower part (for appropriate $r$).
For even $k$, we require $n=4r+2$ to ensure that the difference $\frac{1}{2}(2\sum_{i=2}^k(2i+1)-2)+\frac{n}{2}$ is odd.
For $k=1$, $\overline{T}_{1,n}=T_{1,n}=S_{1,n}=A_6(n)$, i.e., for $n=8+4r$, $r\geq 0$, 
the examples $\overline{T}_{1,n}$ are triangulated tori. 

\begin{thm}
Let $q=6k$ for $k\geq 1$ and let $n=4r$ if $k$ is odd and $n=4r+2$ if $k$ is even, with $r$ an integer
such that $n\geq 3|2\sum_{i=2}^k(2i+1)-2|+2$. Then the series $T_{k,n}$ provides 
infinitely many orientable cyclic triangulations, which are equivelar of type $\{3,q\}=\{3,6k\}$, $k\geq 1$.
\end{thm}

\subsection{The series \mathversion{bold}$U_{k,n}$\mathversion{normal}}

\begin{figure}
\begin{center}
\small
\psfrag{0}{0}
\psfrag{1}{1}
\psfrag{2}{2}
\psfrag{5}{5}
\psfrag{10}{10}
\psfrag{16}{$16\!+\!r$}
\psfrag{22}{$22\!+\!2r=2+\frac{n}{2}$}
\psfrag{[1]}{\scriptsize $\langle 1\rangle$}
\psfrag{[2]}{\scriptsize $\langle 2\rangle$}
\psfrag{[3]}{\scriptsize $\langle 3\rangle$}
\psfrag{[5]}{\scriptsize $\langle 5\rangle$}
\psfrag{[8]}{\scriptsize $\langle 8\rangle$}
\psfrag{[6]}{\scriptsize $\langle 6\!+\!r\rangle$}
\psfrag{[14]}{\scriptsize $\langle 14\!+\!r\rangle$}
\psfrag{[20]}{\scriptsize $\langle 20\!+\!2r\rangle$}
\psfrag{[n/2]}{\scriptsize $\langle \frac{n}{2}\rangle$}
\includegraphics[width=0.5\linewidth]{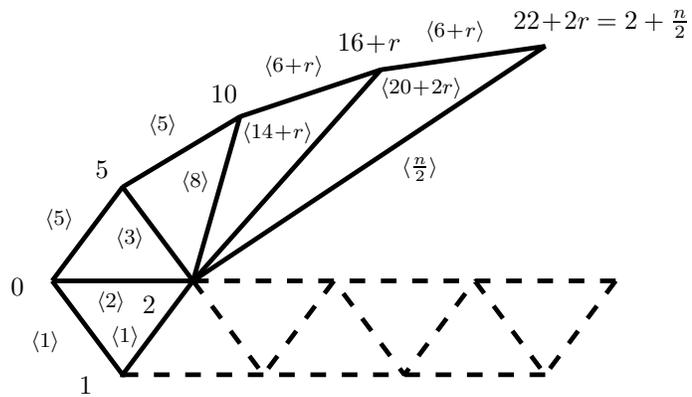}
\end{center}
\caption{The seed of $U_{2,n}$.}
\label{fig:U_2_n}
\end{figure}

\begin{figure}
\begin{center}
\small
\psfrag{0}{0}
\psfrag{1}{1}
\psfrag{-1}{$-1$}
\psfrag{2}{2}
\psfrag{-2}{$-2$}
\psfrag{3}{3}
\psfrag{-3}{$-3$}
\psfrag{5}{5}
\psfrag{-5}{$-5$}
\psfrag{6}{6}
\psfrag{-6}{$-6$}
\psfrag{8}{8}
\psfrag{-8}{$-8$}
\psfrag{10}{10}
\psfrag{-10}{$-10$}
\psfrag{14}{14}
\psfrag{-14}{$-14$}
\psfrag{20}{$20=-20$}
\includegraphics[width=0.25\linewidth]{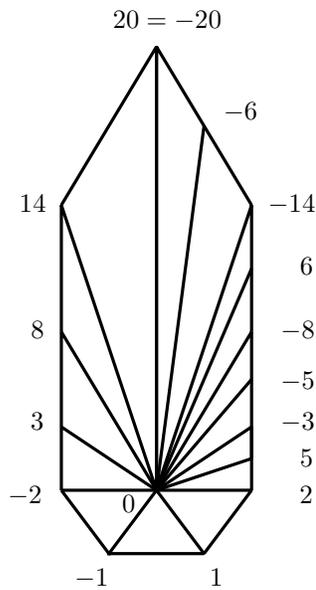}
\end{center}
\caption{The link of $0$ in $U_{2,40}$.}
\label{fig:U_2_40_link}
\end{figure}

Let $q=3+6k$ for $k\geq 1$ and let $n=2(2\sum_{i=2}^{k}(i+3)+2(k+4+r)-2)=2k^2+18k-4+4r$ for $r\geq 0$. 
The series $U_{k,n}$ is defined by the differences $\langle 5\rangle$, $\langle 5\rangle$, $\langle 6\rangle$, $\langle 6\rangle$, \dots, 
$\langle k+3\rangle$, $\langle k+3\rangle$, $\langle k+4+r\rangle$, $\langle k+4+r\rangle$,
with the boundary edge $[2,2\sum_{i=2}^{k}(i+3)+2(k+4+r)]$ generating a short orbit of edges of size $\frac{n}{2}$;
see Figure~\ref{fig:U_2_n} for the seed of $U_{2,n}$ and Figure~\ref{fig:U_2_40_link} for the link of $0$ in $U_{2,40}$.
For $k\geq 2$, the examples $U_{k,n}$ are non-orientable.
For $k=1$, we have $q=9$, $n=2(2(5+r)-2)=16+4r$ for $r\geq 0$, and $U_{1,n}$ is defined by the differences $\langle 5+r\rangle$, $\langle 5+r\rangle$,
with the boundary edge $[2,10+2r]$ generating a short orbit of edges of size $\frac{n}{2}$.
The examples $U_{1,16+4r}$ coincide for $r\geq 0$ with the examples $B_9(16+4r)$ of Section~\ref{subsec:B_9_n}.
The examples $U_{1,n}$ are orientable for $n=16+8s$, $s\geq 0$, and are non-orientable 
for $n=20+8s$, $s\geq 0$.

For $k=1$, $q=9$, $n=16+4r$, $r\geq 0$, let alternatively the series $\overline{U}_{1,n}$ 
be defined by the difference~$\langle \frac{n}{2}\rangle$ and the closing triangles $[2,\frac{n}{4}+1,\frac{n}{2}]$.
The examples of the series $\overline{U}_{1,n}$ are orientable for $n=16+8s$, $s\geq 0$,
and are non-orientable for $n=20+8s$, $s\geq 0$.

\begin{thm}
Let $q=3+6k$ for $k\geq 1$. Then for $k\geq 2$ and $n=2k^2+18k-4+4r$, $r\geq 0$, the series $U_{k,n}$ 
and for $k=1$ and $n=20+8s$, $s\geq 0$, the series $U_{1,n}$ and $\overline{U}_{1,n}$
provide infinitely many non-orientable cyclic triangulations, which are equivelar of type $\{3,q\}=\{3,3+6k\}$.
\end{thm}

\subsection{The series \mathversion{bold}$V_{k,n}$\mathversion{normal}}

Let $q=3+6k$ for $k\geq 1$ and let $n=2(2(\sum_{i=2}^{k+1}(2i+1)+2s)-2)=4k^2+16k-4+8s$ for $s\geq 0$. 
The series $V_{k,n}$ is defined by the odd differences
$\langle 5\rangle$, $\langle 5\rangle$, $\langle 7\rangle$, $\langle 7\rangle$, 
\dots, $\langle 2k+1\rangle$, $\langle 2k+1\rangle$, $\langle 2(k+1)+1+2s\rangle$, $\langle 2(k+1)+1+2s\rangle$,
with the boundary edge $[2,2+\frac{n}{2}]=[2,2k^2+8k+4s]$ generating a short orbit of edges of size $\frac{n}{2}$.
The examples $V_{1,16+8s}$ coincide for $s\geq 0$ with the examples $B_9(16+8s)$ of Section~\ref{subsec:B_9_n}
and with the examples $U_{1,16+8s}$ above.

\begin{thm}
Let $q=3+6k$ for $k\geq 1$. Then for $n=4k^2+16k-4+8s$, $s\geq 0$, the series~$V_{k,n}$ 
provides infinitely many orientable cyclic triangulations, 
which are equivelar of type $\{3,q\}=\{3,3+6k\}$.
\end{thm}

\subsection{The series \mathversion{bold}$W_{k,n}$\mathversion{normal}}

\begin{figure}
\begin{center}
\small
\psfrag{0}{0}
\psfrag{1}{1}
\psfrag{2}{2}
\psfrag{5}{5}
\psfrag{10}{10}
\psfrag{16}{$16\!+\!t$}
\psfrag{22}{$22\!+\!2t=2+\frac{2}{3}n$}
\psfrag{2+2/3n}{$2+\frac{2}{3}n$}
\psfrag{2+1/3n}{$2+\frac{1}{3}n$}
\psfrag{[1]}{\scriptsize $\langle 1\rangle$}
\psfrag{[2]}{\scriptsize $\langle 2\rangle$}
\psfrag{[3]}{\scriptsize $\langle 3\rangle$}
\psfrag{[5]}{\scriptsize $\langle 5\rangle$}
\psfrag{[8]}{\scriptsize $\langle 8\rangle$}
\psfrag{[6]}{\scriptsize $\langle 6\!+\!t\rangle$}
\psfrag{[14]}{\scriptsize $\langle 14\!+\!t\rangle$}
\psfrag{[20]}{\scriptsize $\langle 20\!+\!2t\rangle$}
\psfrag{[2/3n]}{\scriptsize $\langle \frac{2}{3}n\rangle$}
\psfrag{[1/3n]}{\scriptsize $\langle \frac{1}{3}n\rangle$}
\psfrag{[n/2]}{\scriptsize $\langle \frac{n}{2}\rangle$}
\includegraphics[width=0.5\linewidth]{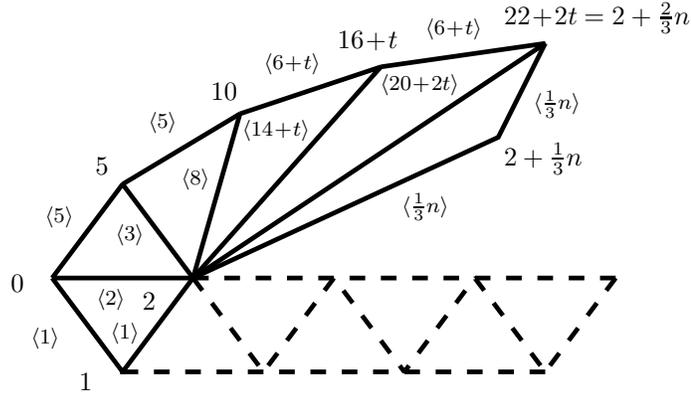}
\end{center}
\caption{The seed of $W_{2,n}$.}
\label{fig:W_2_n}
\end{figure}

\begin{figure}
\begin{center}
\small
\psfrag{0}{0}
\psfrag{1}{1}
\psfrag{-1}{$-1$}
\psfrag{2}{2}
\psfrag{-2}{$-2$}
\psfrag{3}{3}
\psfrag{-3}{$-3$}
\psfrag{5}{5}
\psfrag{-5}{$-5$}
\psfrag{6}{6}
\psfrag{-6}{$-6$}
\psfrag{8}{8}
\psfrag{-8}{$-8$}
\psfrag{10}{10}
\psfrag{-10}{$-10$}
\psfrag{14}{14}
\psfrag{-14}{$-14$}
\psfrag{20}{20}
\psfrag{-20}{$-20$}
\includegraphics[width=0.25\linewidth]{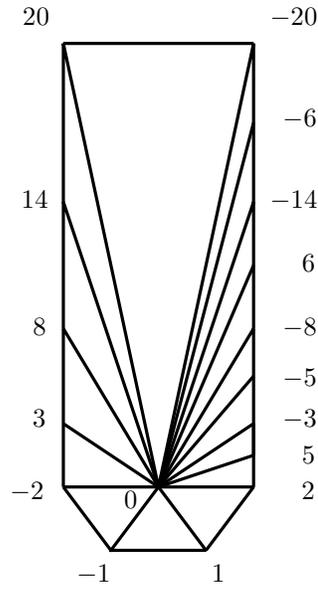}
\end{center}
\caption{The link of $0$ in $W_{2,30}$.}
\label{fig:W_2_30_link}
\end{figure}

\begin{figure}
\begin{center}
\small
\psfrag{0}{0}
\psfrag{1}{1}
\psfrag{2}{2}
\psfrag{3}{3}
\psfrag{4}{4}
\psfrag{5}{5}
\psfrag{6}{6}
\psfrag{7}{7}
\psfrag{8}{8}
\psfrag{9}{9}
\psfrag{10}{10}
\psfrag{11}{11}
\includegraphics[width=0.65\linewidth]{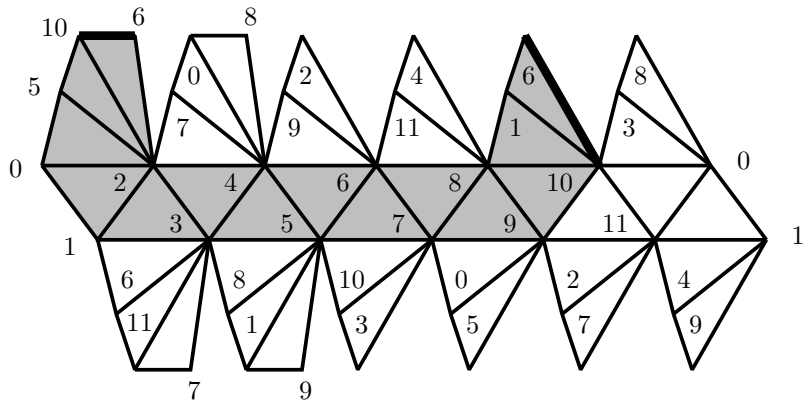}
\end{center}
\caption{The orientable example $W_{1,12}$.}
\label{fig:W_1_12}
\end{figure}

Let the examples $W_{k,n}$ be defined by (the reduced seed given by) the differences 
$\langle 5\rangle$, $\langle 5\rangle$, $\langle 6\rangle$, $\langle 6\rangle$, \dots, $\langle k+3\rangle$, $\langle k+3\rangle$, $\langle k+4+t\rangle$, $\langle k+4+t\rangle$
and the short orbit of the triangle $[2,2+\frac{1}{3}n,2+\frac{2}{3}n]$.
Then $2+\frac{2}{3}n=2\sum_{i=2}^{k+1}(i+3)+2t$, which yields 
$n=\frac{3}{2}(2\sum_{i=2}^{k+1}(i+3)+2t-2)=\frac{3}{2}(k^2+9k+2t-2)$ for $k\geq 1$ and $t\geq 0$; 
see Figure~\ref{fig:W_2_n} for the seed of $W_{2,n}$ and Figure~\ref{fig:W_2_30_link} for the link of $0$
in $W_{2,30}$.

For $t=0$, the examples $W_{k,n}$ coincide with the examples $S_{k+1,n}$, where $n=\frac{3}{2}(k^2+9k-2)$.
Figure~\ref{fig:W_1_12} displays the orientable example $W_{1,12}=S_{2,12}=C_{10}(12)$.
We see in Figure~\ref{fig:W_1_12} that the boundary edges $[2,6]$ and $[6,10]$ of the short (orbit generating) triangle 
$[2,2+\frac{1}{3}n,2+\frac{2}{3}n]=[2,2+\frac{1}{2}(k^2+9k-2),2+(k^2+9k-2)]=[2,6,10]$
in $W_{1,12}$ have \emph{even} differences, so their second copies lie in the upper part. 
Nevertheless, since the short orbit has only $\frac{n}{3}$ triangles, the connecting bands 
between each two copies of the edges are cylinders, as depicted in grey for the band connecting 
the two copies of the edge $[6,10]$ in Figure~\ref{fig:W_1_12}. 

If $k=1$, then $[2,2+\frac{1}{2}(k^2+9k+2t-2),2+(k^2+9k+2t-2)]=[2,6+t,10+2t]$,
with the two boundary edges $[2,6+t]$ and $[6+t,10+2t]$ of difference $4+t$.
The value $4+t$ is even respectively odd, whenever $t$ is even respectively odd. 
Thus, $W_{1,n}$ is orientable for even $t$ and non-orientable for odd $t$.

If $k=2$, then $[2,2+\frac{1}{2}(k^2+9k+2t-2),2+(k^2+9k+2t-2)]=[2,12+t,22+2t]$.
Again, the difference $10+t$ is even respectively odd, whenever $t$ is even respectively odd. 
However, for even $t$, the differences $\langle 5\rangle$, $\langle 5\rangle$, $\langle 6+t\rangle$, $\langle 6+t\rangle$
force the triangulations $W_{2,n}$ to be non-orientable.

\begin{thm}
Let $q=3+6k+1$ for $k\geq 1$. Then for $n=\frac{3}{2}(k^2+9k+2t-2)$, $t\geq 0$, the series $W_{k,n}$ 
provides infinitely many non-orientable cyclic triangulations, which are equivelar 
of type $\{3,q\}=\{3,3+6k+1\}$. For $k\geq 2$, all the examples $W_{k,n}$ are non-orientable.
For $k=1$, $W_{1,n}$ is non-orientable for all odd $t$ and is orientable for all even $t$.
\end{thm}

\subsection{The series \mathversion{bold}$X_{k,n}$\mathversion{normal}}

Let for $k\geq 1$ the examples $X_{k,n}$ be defined by the reduced seed given by the differences 
$\langle 5\rangle$, $\langle 5\rangle$, $\langle 7\rangle$, $\langle 7\rangle$, 
\dots, $\langle 2k+1\rangle$, $\langle 2k+1\rangle$, $\langle 2k+3+2t\rangle$, $\langle 2k+3+2t\rangle$
and the short orbit of the triangle $[2,2+\frac{1}{3}n,2+\frac{2}{3}n]$.
Then $2+\frac{2}{3}n=2\sum_{i=1}^{k}(2i+3)+4t$, which yields 
$n=\frac{3}{2}(2\sum_{i=1}^{k}(2i+3)+4t-2)=3(k^2+4k+2t-1)$ for $k\geq 1$ and $t\geq 0$.
For $k$ odd and $t\geq 0$, the examples $X_{k,n}$ are orientable;
for $k$ odd and $t=0$, the examples $X_{k,n}$ coincide with the examples $T_{k+1,n}$, where $n=3(k^2+4k-1)$.
For $k$ even, $n=3(k^2+4k+2t-1)$ is odd for all $t$, implying that $X_{k,n}$ is non-orientable.

\begin{thm}
Let $q=9+12\ell+1$ for $k=2\ell+1$, $\ell\geq 0$. Then for $n=3(k^2+4k+2t-1)$, $t\geq 0$, the series $X_{k,n}$ 
provides infinitely many orientable cyclic triangulations, which are equivelar 
of type $\{3,q\}=\{3,9+12\ell+1\}$.
\end{thm}

\subsection{The series \mathversion{bold}$\overline{X}_{k,n}$\mathversion{normal}}
\label{subsec:XX}

\begin{figure}
\begin{center}
\small
\psfrag{0}{0}
\psfrag{1}{1}
\psfrag{2}{2}
\psfrag{4}{4}
\psfrag{5}{5}
\psfrag{10}{10}
\psfrag{17+2r}{$17\!+\!2t$}
\psfrag{24+4r}{$24\!+\!4t$}
\psfrag{44+8r}{$44\!+\!8t$}
\psfrag{23+4r}{$23\!+\!4t$}
\psfrag{[1]}{\scriptsize $\langle 1\rangle$}
\psfrag{[2]}{\scriptsize $\langle 2\rangle$}
\psfrag{[3]}{\scriptsize $\langle 3\rangle$}
\psfrag{[5]}{\scriptsize $\langle 5\rangle$}
\psfrag{[8]}{\scriptsize $\langle 8\rangle$}
\psfrag{[7+2r]}{\scriptsize $\langle 7\!+\!2t\rangle$}
\psfrag{[15+2r]}{\scriptsize $\langle 15\!+\!2t\rangle$}
\psfrag{[22+4r]}{\scriptsize $\langle 22\!+\!4t\rangle$}
\psfrag{[21+4r]}{\scriptsize $\langle 21\!+\!4t\rangle$}
\psfrag{[42+8r]}{\scriptsize $\langle 42\!+\!8t\rangle$}
\psfrag{[20+4r]}{\scriptsize $\langle 20\!+\!4t\rangle$}
\psfrag{[-20-4r]}{\scriptsize $\langle -20\!-\!4t\rangle$}
\includegraphics[width=0.65\linewidth]{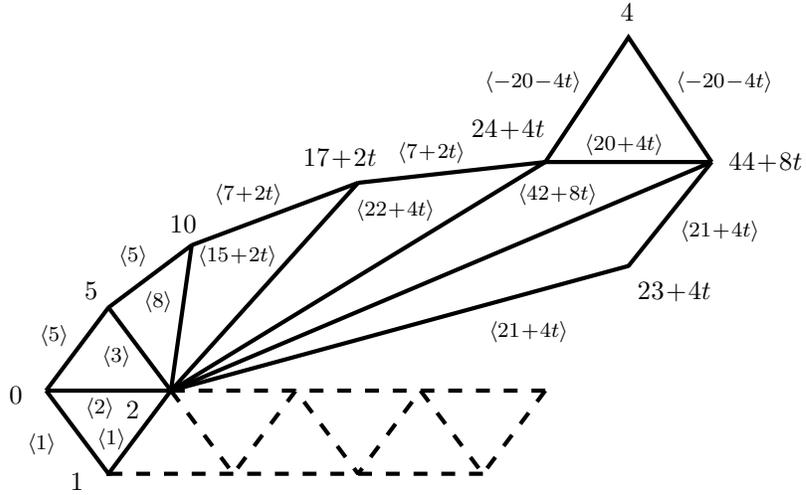}
\end{center}
\caption{The seed of $\overline{X}_{3,n}$.}
\label{fig:XX_3_n}
\end{figure}

Let $n=3(2\sum_{i=2}^{k}(2i+1)+4t-4)$ for $k\geq 2$, $t\geq 0$. The examples $\overline{X}_{k,n}$ are defined 
by the differences 
$\langle 5\rangle$, $\langle 5\rangle$, $\langle 7\rangle$, $\langle 7\rangle$, 
\dots, $\langle 2k-1\rangle$, $\langle 2k-1\rangle$, $\langle 2k+1+2t\rangle$, $\langle 2k+1+2t\rangle$
and the three triangles 
$[2,2\sum_{i=2}^{k}(2i+1)+4t,4\sum_{i=2}^{k}(2i+1)+4t-4]$,
$[2,2\sum_{i=2}^{k}(2i+1)+4t-1,4\sum_{i=2}^{k}(2i+1)+4t-4]$,
$[4,2\sum_{i=2}^{k}(2i+1)+4t,4\sum_{i=2}^{k}(2i+1)+4t-4]$,
of which the last triangle generates a short orbit;
see Figure~\ref{fig:XX_3_n} for the seed of $\overline{X}_{3,n}$.

\begin{thm}
Let $q=3+6k+1$ for $k\geq 2$. Then for $n=3(2\sum_{i=2}^{k}(2i+1)+4t-4)$, $t\geq 0$, the series $\overline{X}_{k,n}$ 
provides infinitely many orientable cyclic triangulations, which are equivelar 
of type $\{3,q\}=\{3,3+6k+1\}$.
\end{thm}

\subsection{The series \mathversion{bold}$Y_{k,n}$\mathversion{normal}}

\begin{figure}
\begin{center}
\small
\psfrag{0}{0}
\psfrag{1}{1}
\psfrag{2}{2}
\psfrag{5}{5}
\psfrag{10}{10}
\psfrag{16}{$16\!+\!t$}
\psfrag{22}{$22\!+\!2t$}
\psfrag{2+2/3n}{$22+2t+\frac{n}{2}=2+\frac{2}{3}n$}
\psfrag{2+1/3n}{$2+\frac{1}{3}n$}
\psfrag{[1]}{\scriptsize $\langle 1\rangle$}
\psfrag{[2]}{\scriptsize $\langle 2\rangle$}
\psfrag{[3]}{\scriptsize $\langle 3\rangle$}
\psfrag{[5]}{\scriptsize $\langle 5\rangle$}
\psfrag{[8]}{\scriptsize $\langle 8\rangle$}
\psfrag{[6]}{\scriptsize $\langle 6\!+\!t\rangle$}
\psfrag{[14]}{\scriptsize $\langle 14\!+\!t\rangle$}
\psfrag{[20]}{\scriptsize $\langle 20\!+\!2t\rangle$}
\psfrag{[2/3n]}{\scriptsize $\langle \frac{2}{3}n\rangle$}
\psfrag{[1/3n]}{\scriptsize $\langle \frac{1}{3}n\rangle$}
\psfrag{[n/2]}{\scriptsize $\langle \frac{n}{2}\rangle$}
\includegraphics[width=0.625\linewidth]{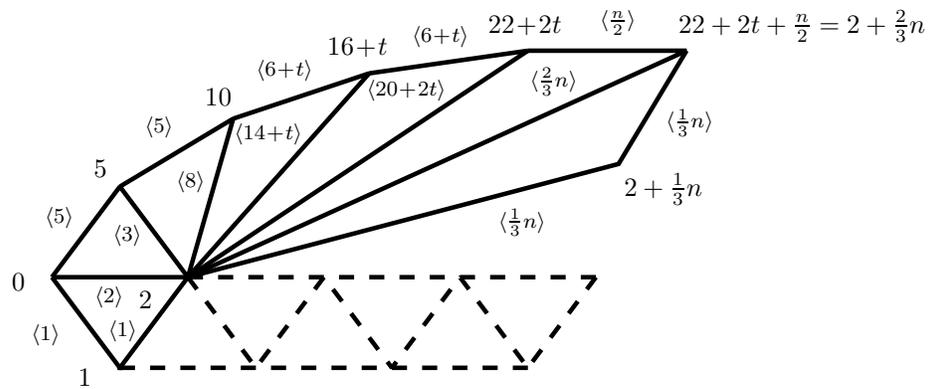}
\end{center}
\caption{The seed of $Y_{3,n}$.}
\label{fig:Y_3_n}
\end{figure}

\begin{figure}
\begin{center}
\small
\psfrag{0}{0}
\psfrag{1}{1}
\psfrag{-1}{$-1$}
\psfrag{2}{2}
\psfrag{-2}{$-2$}
\psfrag{3}{3}
\psfrag{-3}{$-3$}
\psfrag{5}{5}
\psfrag{-5}{$-5$}
\psfrag{6}{6}
\psfrag{-6}{$-6$}
\psfrag{8}{8}
\psfrag{-8}{$-8$}
\psfrag{10}{10}
\psfrag{-10}{$-10$}
\psfrag{14}{14}
\psfrag{-14}{$-14$}
\psfrag{20}{20}
\psfrag{-20}{$-20$}
\psfrag{40}{\scriptsize $40=-40$}
\psfrag{60}{\scriptsize $60=-60$}
\psfrag{80}{\scriptsize $80=-80$}
\includegraphics[width=0.25\linewidth]{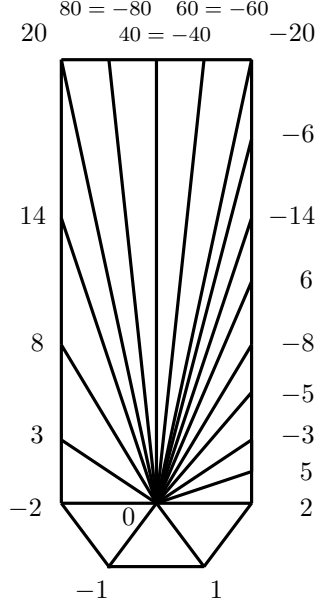}
\end{center}
\caption{The link of $0$ in $Y_{3,120}$.}
\label{fig:Y_3_120_link}
\end{figure}

Let for $k\geq 2$ and $t\geq 0$ the examples $Y_{k,n}$ be defined by the differences 
$\langle 5\rangle$, $\langle 5\rangle$, $\langle 6\rangle$, $\langle 6\rangle$, 
\dots, $\langle k+2\rangle$, $\langle k+2\rangle$, $\langle k+3+t\rangle$, $\langle k+3+t\rangle$, $\langle \frac{n}{2}\rangle$
and the short orbit $[2,2+\frac{1}{3}n,2+\frac{2}{3}n]$,
where $2+\frac{2}{3}n=2\sum_{i=2}^k(i+3)+2t+\frac{n}{2}$, that is, $n=6(2\sum_{i=2}^k(i+3)+2t-2)=6(k^2+7k+2t-10)$.
The seed of $Y_{3,n}$ is displayed in Figure~\ref{fig:Y_3_n}, for the link of $0$ in $Y_{3,120}$ see Figure~\ref{fig:Y_3_120_link}.

\begin{thm}
Let $q=6k+1$ for $k\geq 2$ and let $n=6(k^2+7k+2t-10)$ for $t\geq 0$. 
Then for $k\geq 4$ the series $Y_{k,n}$ provides infinitely many non-orientable cyclic triangulations, 
which are equivelar of type $\{3,q\}=\{3,6k+1\}$.
For $k=2$, the series $Y_{k,n}$ provides non-orientable examples when $t\geq 0$ is odd,
whereas for $k=3$, the series $Y_{k,n}$ provides non-orientable examples when $t\geq 0$ is even.
\end{thm}

\subsection{The series \mathversion{bold}$Z_{k,n}$\mathversion{normal}}

Let for $k\geq 2$ and $t\geq 0$ the examples $Z_{k,n}$ be defined by the differences 
$\langle 5\rangle$, $\langle 5\rangle$, $\langle 7\rangle$, $\langle 7\rangle$, \dots, 
$\langle 2k-1\rangle$, $\langle 2k-1\rangle$, $\langle 2k+1+2t\rangle$, $\langle 2k+1+2t\rangle$, $\langle \frac{n}{2}\rangle$
and the short orbit $[2,2+\frac{1}{3}n,2+\frac{2}{3}n]$.
Then $2+\frac{2}{3}n=2\sum_{i=2}^k(2i+1)+4t+\frac{n}{2}$, 
from which it follows that $n=6(2k^2+4k+4t-8)$.

\begin{thm}
Let $q=6k+1$ for $k\geq 2$ and let $n=6(2k^2+4k+4t-8)$ for $t\geq 0$. 
Then for $k\geq 2$ the series $Z_{k,n}$ provides infinitely many orientable cyclic triangulations, 
which are equivelar of type $\{3,q\}=\{3,6k+1\}$.
\end{thm}

\subsection{The case \mathversion{bold}$q=7$\mathversion{normal}}

\begin{figure}
\begin{center}
\small
\psfrag{0}{0}
\psfrag{+-k}{$\pm j(+\frac{n}{2})$}
\psfrag{n/2}{$\frac{n}{2}$}
\psfrag{+-n/3}{$\pm\frac{n}{3}$}
\psfrag{-+n/3=+-2k+n/2}{$\mp\frac{n}{3}=\pm 2j+\frac{n}{2}(+n)$}
\psfrag{+-k+n/2}{$\pm j+\frac{n}{2}(+\frac{n}{2})$}
\psfrag{[+-k]}{\scriptsize $\langle\pm j(+\frac{n}{2})\rangle$}
\psfrag{[+-n/3]}{\scriptsize $\langle\pm\frac{n}{3}\rangle$}
\psfrag{[-+n/3]}{\scriptsize $\langle\mp\frac{n}{3}\rangle$}
\psfrag{[n/2]}{\scriptsize $\langle\frac{n}{2}\rangle$}
\includegraphics[width=0.85\linewidth]{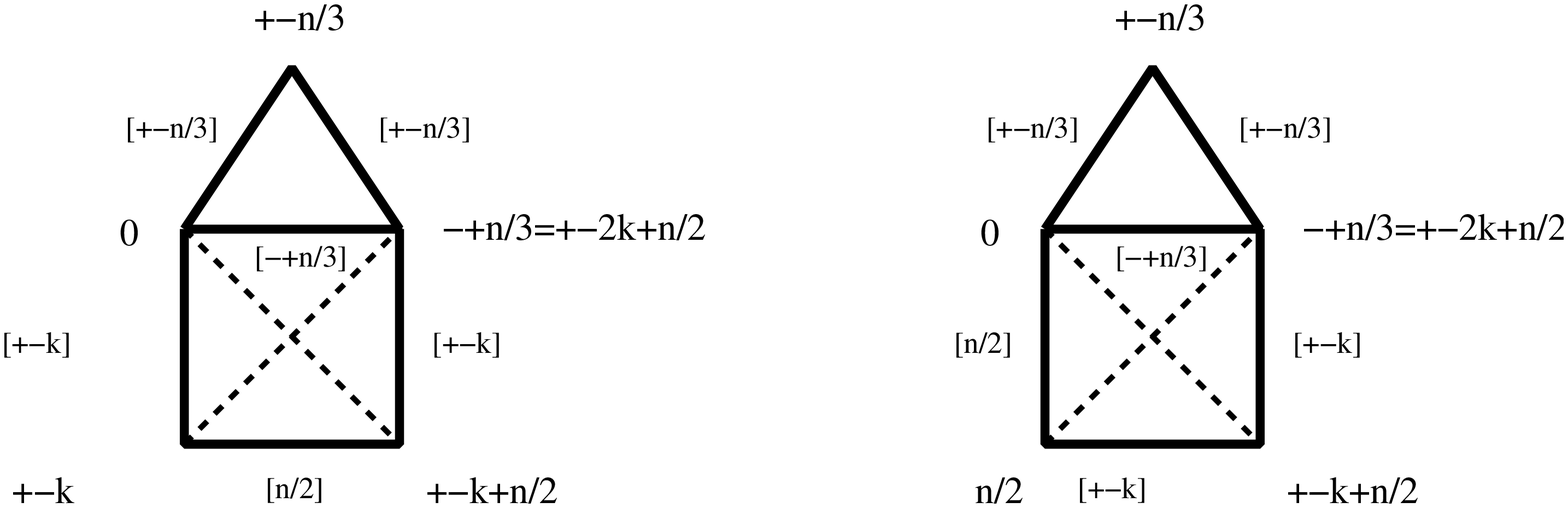}
\end{center}
\caption{The case $q=7$: possible configurations.}
\label{fig:q7}
\end{figure}

In the case of $q=7$, there are three orbits, two of size $n$ and one short orbit of size $\frac{n}{2}$.
In particular, any fundamental domain consists of three triangles. By a cyclic shift,
we can choose the short triangle to be $[0,\pm\frac{n}{3},\mp\frac{n}{3}]$.
The other two triangles of the fundamental domain form a square and are glued,
say, via the edge $[0,\mp\frac{n}{3}]$ to the short triangle. There are three
remaining boundary edges of the square, so one of these edges must be of length $\frac{n}{2}$,
while the other two edges both have length $\pm j(+\frac{n}{2})$. In Figure~\ref{fig:q7}, we see all possible
configurations, with the dashed lines indicating the two choices for the diagonals of the square.
Either the two edges of length $\pm j(+\frac{n}{2})$ are consecutive (in the right of Figure~\ref{fig:q7})
or they are separated by the edge of length $\frac{n}{2}$ (in the left of Figure~\ref{fig:q7}).
However, in both situations, we have the closing condition $\mp\frac{n}{3}=\pm 2j+\frac{n}{2}$,
which either gives $n=12j$ or $n=\frac{12}{5}j$. 
In the first case, all occurring differences and vertices in the diagram are multiples of $j$, 
which implies that the resulting surface is not connected for $j>1$. Hence, $j=1$, and therefore $n=12$.
Analogously, it follows in the second case that $j=5$, and again $n=12$.

It is then easy to see, that there are exactly two cyclic triangulations with $q=7$ and $n=12$,
both of the orientable surface of genus $g=2$, the one with generating triangles $[0,1,2]$, $[0,2,6]$, $[0,4,8]$
and the other with generating triangles $[0,1,5]$, $[0,1,6]$, $[0,4,8]$.

\begin{thm}
For $q=7$, there are exactly two cyclic $7$-equivelar triangulations, both of the
orientable surface of genus $g=2$ with $n=12$ vertices.
\end{thm}

\subsection{Examples with non-pathlike seeds}

\begin{figure}
\begin{center}
\small
\psfrag{0}{0}
\psfrag{1}{1}
\psfrag{2}{2}
\psfrag{5}{5}
\psfrag{n+5/2}{$\frac{n+5}{2}$}
\psfrag{n+7/2}{$\frac{n+7}{2}$}
\psfrag{[1]}{\scriptsize $\langle 1\rangle$}
\psfrag{[2]}{\scriptsize $\langle 2\rangle$}
\psfrag{[3]}{\scriptsize $\langle 3\rangle$}
\psfrag{[5]}{\scriptsize $\langle 5\rangle$}
\psfrag{[n+3/2]}{\scriptsize $\langle \frac{n+3}{2}\rangle$}
\psfrag{[n+5/2]}{\scriptsize $\langle \frac{n+5}{2}\rangle$}
\includegraphics[width=0.6\linewidth]{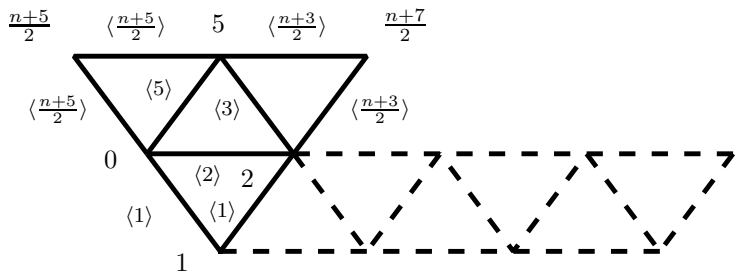}
\end{center}
\caption{The seed of $E_{18}(n)$.}
\label{fig:E_18_n}
\end{figure}

Not all members of the $[0,1,2]$-family have a pathlike seed of triangles. 
The series $E_{18}(n)$ of Section~\ref{subsec:E_12_n}, which is defined for $n=17+2m$, $m\geq 0$, 
and the series $F_{18}(n)$ of Section~\ref{subsec:F_12_n}, which is defined for $n=14$ and $n=18+2m$, $m\geq 0$, 
provide such examples that do not belong to the family $\tilde{G}_{d_1,d_1,d_2,d_2,\dots,d_{k-1},d_{k-1}}$. 
The seed of $E_{18}(n)$ is displayed in Figure~\ref{fig:E_18_n}. Furthermore, the examples 
of the two-parameter series $\overline{X}_{k,n}$ of Section~\ref{subsec:XX} have non-pathlike seeds.

\section{Cyclic triangulations beyond the \mathversion{bold}$[0,1,2]$-triangulations\mathversion{normal}-family}
\label{sec:beyond}

The examples of the orientable cyclic series $G_{18}(n)$ of $\{3,18\}$-equivelar triangulations from Section~\ref{subsec:G_18_n_H_18_n}, 
defined for $n=19+2m$, $m\geq 0$, do not belong to the $[0,1,2]$-family. 
Also, the examples of the two-parameter generalized Ringel series $R_{k,n}$ 
of cyclic triangulations \cite{LutzSulankeTiwariUpadhyay2010pre} with $k\geq 0$, $n\geq 7+12k$ 
do not belong to the $[0,1,2]$-family. Nevertheless, the examples $G_{18}(n)$ (and also the examples $R_{k,n}$) 
contain the orbit of the edge $[0,1]$. In the case of the series $G_{18}(n)$,
we used the cycle $0$--$1$--$2$-- $\cdots$ --$(n-1)$--$0$ as a ``horizontal cycle'' for a symmetric display 
of $G_{18}(19)$ in Figure~\ref{fig:G_18_19}. The image of $G_{18}(19)$ can be cut into 
fundamental domains by the cyclic translates of the edge $[0,9]$.
Any such fundamental domain is a triangulated disk with boundary edges (in clockwise order) 
of differences 
$\langle 9\rangle$, $\langle -6\rangle$, $\langle 7\rangle$, $\langle -9\rangle$,
$\langle -8\rangle$, $\langle -7\rangle$, $\langle 8\rangle$, and $\langle 6\rangle$.

In fact, the fundamental domain of any standard cyclic triangulation with $q=3k$, $k\geq 1$, 
is a triangulated disk for which the boundary edges can be grouped in pairs according
to their differences if the number of boundary edges is even. If the number of boundary edges is odd,
then there is one boundary edge of difference $\frac{n}{2}$, whereas all other boundary edges
can be grouped in pairs according to their differences. In the non-standard case, $q=3k+1$, $k\geq 1$,
the description of the fundamental domain has to be modified to accommodate a short orbit of triangles
of size~$\frac{n}{3}$.

\section{Equivelar tessellations of surfaces with cyclic symmetry}
\label{sec:tessellations}

Let a \emph{map} on a surface $M$ be a decomposition of $M$ into a finite cell complex 
and let $G$ be the $1$-skeleton of the map on $M$. The graph $G$ of the map may have 
multiple edges, loops, vertices of degree $2$, or even vertices of degree $1$. 
For example, the embedding of a tree with $n$ vertices and $n-1$ edges on $S^2$ 
decomposes the $2$-sphere into one polygon with $2n-2$ edges, which are identified pairwise. 
(Sometimes the graphs of maps are required to be connected finite
simple graphs, sometimes multiple edges are allowed but no loops,
and vertices are often required to have at least degree $3$;
see \cite{BrehmWills1993}, \cite{CoxeterMoser1957}, \cite{Ringel1974},~\cite{Vince2004}.)

A map is \emph{equivelar of type $\{p,q\}$} if $M$ is decomposed into
$p$-gons only with every vertex having degree $q$; cf.\
\cite{McMullenSchulzWills1982}, \cite{McMullenSchulzWills1983}. 
A map is \emph{polyhedral} if the intersection of any two of its polygons is either empty, a common vertex, 
or a common edge; see the surveys \ \cite{BrehmSchulte1997}, \cite{BrehmWills1993}.
An \emph{equivelar polyhedral map} is a map which is both equivelar and polyhedral. 

\begin{figure}
\begin{center}
\small
\psfrag{0}{0}
\psfrag{1}{1}
\psfrag{2}{2}
\psfrag{3}{3}
\psfrag{4}{4}
\psfrag{5}{5}
\psfrag{6}{6}
\psfrag{7}{7}
\psfrag{8}{8}
\psfrag{9}{9}
\psfrag{10}{10}
\psfrag{11}{11}
\psfrag{12}{12}
\psfrag{13}{13}
\psfrag{14}{14}
\psfrag{15}{15}
\includegraphics[width=0.72\linewidth]{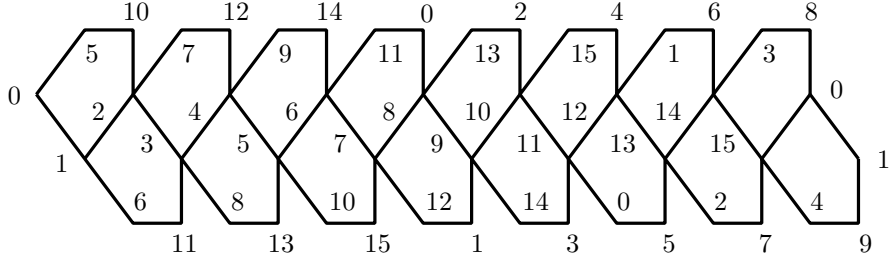}
\end{center}
\caption{The $\{5,5\}$-equivelar example $B_{\{5,5\}}(16)$.}
\label{fig:B_5_5_16}
\end{figure}

\begin{thm}
The deletion of orbits of interior seed edges from any standard cyclic example 
of the family $\tilde{G}_{d_1,d_1,d_2,d_2,\dots,d_{k-1},d_{k-1}}$ of Section~\ref{sec:G}
yields a \emph{regular} cell decomposition (i.e., a cell decomposition without identifications 
on the boundaries of the cells) of the respective surface or pinched surface. 
If the original triangulation of the (pinched) surface is $q$-equivelar,
then the deletion of all orbits of interior edges gives a tessellated surface,
which is $\{\frac{q}{3}+2,\frac{q}{3}+2\}$-equivelar.
\end{thm}

\textbf{Proof:} 
The seed of a standard example of the $\tilde{G}_{d_1,d_1,d_2,d_2,\dots,d_{k-1},d_{k-1}}$-family is a subcomplex 
of the star of the vertex $2$ and thus contains every vertex on its boundary exactly once, from which the regularity 
of the resulting cell complex follows.

If the example of the $\tilde{G}_{d_1,d_1,d_2,d_2,\dots,d_{k-1},d_{k-1}}$-family is $q$-equivelar,
then exactly $\frac{q}{3}-1$ orbits of interior edges are deleted, so the boundary of a respective seed
forms a $p'$-gon with $p'=\frac{q}{3}+2$. The deletion of $\frac{q}{3}-1$ interior edges reduces the degree
of every vertex by $2(\frac{q}{3}-1)$, hence giving $q'=q-2(\frac{q}{3}-1)=\frac{q}{3}+2$.
The resulting tessellated surface therefore is $\{p',q'\}$-equivelar with 
$\{p',q'\}=\{\frac{q}{3}+2,\frac{q}{3}+2\}$.~\hfill$\Box$

\medskip
\smallskip

If we delete the orbits of the edges $[0,2]$ and $[2,5]$ from example $B_{9}(16)$ (Figure~\ref{fig:B_9_16}) of Section~\ref{subsec:B_9_n},
we obtain an orientable $\{5,5\}$-equivelar polyhedral map $B_{\{5,5\}}(16)$ with $16$~vertices as displayed in Figure~\ref{fig:B_5_5_16}.
A first description of a $\{5,5\}$-equivelar polyhedral map with $16$~vertices is due to Brehm~\cite{Brehm1990}.
Although his description differs from ours, the respective maps are combinatorially isomorphic.
Brehm uses vertices $1$, $2$, \dots, $16$ and differences $\langle 8\rangle$, $\langle -3\rangle$, $\langle -3\rangle$,
whereas we use differences $\langle 5\rangle$,~$\langle 5\rangle$, $\langle -8\rangle$,
where, for brevity, we omitted to note the last difference $\langle -8\rangle$ in Section~\ref{sec:012_series}.
The combinatorially isomorphism that maps our labeling to Brehm's labeling is given by
$0\mapsto 1$, $1\mapsto 4$, $2\mapsto 7$, $3\mapsto 10$, $4\mapsto 13$, $5\mapsto 16$, $6\mapsto 3$, $7\mapsto 6$,
$8\mapsto 9$, $9\mapsto 12$, $10\mapsto 15$, $11\mapsto 2$, $12\mapsto 5$, $13\mapsto 8$, $14\mapsto 11$, $15\mapsto 14$.
The Brehm map has $40$ edges, which is the minimal number of edges that a polyhedral map on the orientable surface 
of Euler characteristic~$-8$ can have~\cite{Brehm1990}. It was proved in \cite{BrehmDattaNilakantan2002}
that the Brehm map is the unique $\{5,5\}$-equivelar polyhedral map with $16$~vertices on a (orientable or non-orientable) 
surface of Euler characteristic~$-8$, which implies that at least $41$ edges are needed for a polyhedral map 
on the non-orientable surface of Euler characteristic $-8$.

The Brehm example was generalized to an infinite series $BDN_{\{5,5\}}(n)$ of $\{5,5\}$-equivelar polyhedral maps 
by Brehm, Datta, and Nilakantan \cite{BrehmDattaNilakantan2002}. The series $BDN_{\{5,5\}}(n)$ is defined 
by the differences $\langle 8+2r\rangle$, $\langle -(3+r)\rangle$, and $\langle -(3+r)\rangle$
and has $n=16+4r$ vertices, where $r\geq 0$. A~further generalization to a two parameter series $D_{\{2k+1,2k+1\}}(n)$
of $\{2k+1,2k+1\}$-equivelar surfaces for $k\geq 2$ is given in~\cite{Datta2005}. The examples of
this series have $n=2\cdot 3^k-2+4r$ vertices and are defined by the differences $\langle 3^k-1+2r\rangle$, 
$\langle -(3^{k-1}+r)\rangle$, $\langle -(3^{k-1}+r)\rangle$, $\langle -3^{k-2}\rangle$, $\langle -3^{k-2}\rangle$, \dots 
$\langle -9\rangle$, $\langle -9\rangle$, $\langle -3\rangle$, $\langle -3\rangle$. For $k=2$, $D_{\{5,5\}}(n)=BDN_{\{5,5\}}(n)$.

\begin{cor}
The $\{2k+1,2k+1\}$-equivelar polyhedral maps $D_{\{2k+1,2k+1\}}(n)$ with $k\geq 2$ and $n=2\cdot 3^k-2+4r$, $r\geq 0$,
are orientable for even $r$ and are non-orientable for odd $r$.
\end{cor}

Besides the $\{5,5\}$-equivelar example, Brehm presented in \cite{Brehm1990} an edge-minimal $\{6,6\}$-equivelar polyhedral map
of the orientable surface of Euler characteristic $-26$ with $26$ vertices,
defined by the differences $\langle -9\rangle$,  $\langle -9\rangle$, $\langle -3\rangle$, $\langle -3\rangle$.
A generalization of this example to a two-parameter series $D_{\{2k,2k\}}(n)$
of $\{2k,2k\}$-equivelar surfaces with $k\geq 3$ and $n=3^k-1+2r$ vertices, $r\geq 0$,
can be found in~\cite{Datta2005}. The examples of the series 
are defined by the differences
$\langle -(3^{k-1}+r)\rangle$, $\langle -(3^{k-1}+r)\rangle$, $\langle -3^{k-2}\rangle$, $\langle -3^{k-2}\rangle$, \dots 
$\langle -9\rangle$, $\langle -9\rangle$, $\langle -3\rangle$, $\langle -3\rangle$. 

\begin{cor}
The $\{2k,2k\}$-equivelar polyhedral maps $D_{\{2k,2k\}}(n)$ with $k\geq 3$ and $n=3^k-1+2r$, $r\geq 0$,
are orientable for even $r$ and are non-orientable for odd $r$.
\end{cor}

For $5\leq k\leq 12$, Jamet lists in his Diploma thesis \cite{Jamet2001} (advised by Ulrich Brehm)
those numbers of vertices $n$ for which $\{k,k\}$-equivelar \emph{polyhedral} maps exist that have
vertex-transitive and facet-transitive cyclic symmetry. A necessary condition for examples in this special class
of cyclic maps is that the boundary edges of a fundamental domain are grouped in pairs of \emph{successive} edges 
of equal distances if the number of boundary edges is even.
If the number of boundary edges is odd, then there is one additional boundary edge of difference $\frac{n}{2}$.
For even $6\leq k\leq 100$, Jamet also provides examples of $\{k,k\}$-equivelar polyhedral 
maps with an edge-transitive symmetry group that contains a vertex-transitive cyclic subgroup. 
For all the examples of cyclic $\{k,k\}$-equivelar \emph{polyhedral} maps, an arbitrary triangulation of one of the $k$-gons
yields a fundamental domain for a cyclic triangulation of the respective surface.

If the seed of a standard example of the $\tilde{G}_{d_1,d_1,d_2,d_2,\dots,d_{k-1},d_{k-1}}$-family consists of
an even number of triangles, then by deleting every second interior edge (where we start with the first interior edge)
we obtain a quadrangulation of the respective surface or pinched surface. E.g., in the case of the series $S_{4,n}$ 
from Section~\ref{subsec:S_k_n}, we remove the orbits of the edges $[0,2]$, $[2,10]$, $[2,22]$, and $[2,36]$ from~$S_{4,n}$. 
Similarly, if the seed contains $t$ triangles and $t$ is divisible by $s$, we can merge successive groups of $s$ triangles, 
which gives a regular cell decomposition into $(s+2)$-gons. 

In general, the resulting regular cell decompositions will not be \emph{strongly regular}  
(i.e., the intersection of any two polygons is either empty, a common vertex, or a common edge,
which precisely is the requirement on a cell decomposition to be a polyhedral map).
E.g., if we delete \emph{all} interior seed edges and their cyclic translates from the examples 
of the series $S_{4,n}$, then the diagonals $[2,5]$ and $[16,19]$ of the resulting $10$-gon
are cyclic translates of each other and prevent the tessellation from being strongly regular.
Nevertheless, the series $D_{\{2k+1,2k+1\}}(n)$ and $D_{\{2k,2k\}}(n)$ as well as the
examples in \cite{Jamet2001} show that it is possible 
to obtain strongly regular tessellations by deleting orbits of interior seed edges
from members of the $\tilde{G}_{d_1,d_1,d_2,d_2,\dots,d_{k-1},d_{k-1}}$-family. 

\bibliography{.}

\bigskip
\bigskip
\medskip

\small

\noindent
Frank H. Lutz\\
Institut f\"ur Mathematik\\
Technische Universit\"at Berlin\\
Stra\ss e des 17.\ Juni 136\\
10623 Berlin, Germany\\
{\tt lutz@math.tu-berlin.de}

\end{document}